\documentclass[12pt]{article}
 \usepackage{graphicx}
 \usepackage{amssymb}
 \usepackage{epstopdf}
 \newtheorem{thm}{Theorem}[section]
  \newtheorem{lem}[thm]{Lemma}
  \newtheorem{cor}[thm]{Corollary}
  \newtheorem{prop}[thm]{Proposition}

  \newtheorem{rem}[thm]{Remark}

  \newtheorem{ques}[thm]{Question}

 \def \Z{\mathbb Z}
 
 \def \HH {\mathrm{H}}

 \def \sT {\mathcal T}
 \def \R{\mathbb R}
 \def\sE{{\mathcal E}}
 \def \sF {{\mathcal F}}
 \def \sS{\mathcal S}
 
 \def\sK{\mathcal K}
 \def\sB{{\mathcal B}}
 \def\sM{{\mathcal M}}
 \def\sQ{{\mathcal Q}}
 
 \def\sR{{\mathcal R}}
 \def\DE{{\mathcal{E}_{\mathrm{dir}}}}
 \def\proj{\pi}
 
 \def\g{\mathfrak g}
 \def\f{\mathfrak f}
 \def\h{\mathfrak h}
 \def\r{\mathfrak r}
 \def\qed{\hfill\framebox(5,5){}}
  \def\proof{\noindent{\it Proof.\ }}

 \title{A family of pseudo-Anosov braids with small dilatation}
\author{
Eriko Hironaka and Eiko Kin
 }
 
\begin{document}
 \maketitle
\abstract{This paper concerns a family of pseudo-Anosov braids with dilatations
arbitrarily close to one. The associated graph maps and train tracks 
have stable ``star-like" shapes, and the characteristic polynomials 
of their transition matrices form Salem-Boyd sequences.
These examples show that the logarithms of least dilatations of 
pseudo-Anosov braids on $2g+1$ strands
are bounded above by $\log(2 + \sqrt{3})/g$.   It follows that the
asymptotic behavior of least dilatations of pseudo-Anosov,
hyperelliptic surface homeomorphisms
is identical to that found by Penner for general surface homeomorphisms.
}

\section{Introduction}

A braid $\beta$ on $s$ strands is {\it pseudo-Anosov} if its associated mapping
class $\phi$ on $S^2$ with $s+1$ marked points is pseudo-Anosov.  In this case
its dilatation $\lambda(\beta)$ is defined to be the dilatation $\lambda(\phi)$ of $\phi$.    
For fixed $s$, the set of dilatations of pseudo-Anosov braids on $s$ strands
consists of real algebraic integers greater than one with bounded degree, and hence 
has a well-defined minimum.
The general problem of finding the least dilatation of a pseudo-Anosov braids on
$s$ strands is open, and exact results are known only for $s \leq 5$.
Let $\sigma_1,\dots,\sigma_{s-1}$ be the standard braid generators for
braids with $s$ strands.
For $s=3$,  Matsuoka \cite{Ma86} showed that the smallest dilatation is realized by
$\sigma_1\sigma_2^{-1}$.  Ko, Los and Song  \cite{SoKoLos02} showed that for 
$s=4$ the smallest dilatation is realized by
$\sigma_1\sigma_2\sigma_3^{-1}$.  Recently Ham and Song have announced
a revised proof for the case $s=4$, and a proof that 
\begin{eqnarray}\label{HamSong-braid}
\sigma_1\sigma_2\sigma_3\sigma_4\sigma_1\sigma_2
\end{eqnarray}
has smallest dilatation for $s=5$ \cite{HamSong05}.

In this paper, we study a family of generalizations of these examples to arbitrary 
numbers of strands.   Let $\sB(D,s)$ denote the braid group on $s$ strands, 
where $D$ denotes the 2-dimensional disk.
First consider the braids $\beta_{m,n}$ in $\sB(D,m+n+1)$ given by
$$
\beta_{m,n} = \sigma_1\cdots\sigma_m\sigma_{m+1}^{-1}\cdots \sigma_{m+n}^{-1}.
$$
Matsuoka's example appears as $\beta_{1,1}$, and Ko, Los and Song's example as
$\beta_{2,1}$.   For any $m,n \ge 1$ the braids $\beta_{m,n}$
are pseudo-Anosov (see Theorem~\ref{Bmn-classification-thm}).
The dilatations of $\beta_{m,m}$ coincide with those found by Brinkman \cite{Brinkmann04}
(see Section~\ref{fiberedQ-section}),  who also shows that the dilatations
arising in this family can be made arbitrarily close to 1. 
\begin{figure}[htbp]
\begin{center}
\includegraphics[height=1.5in]{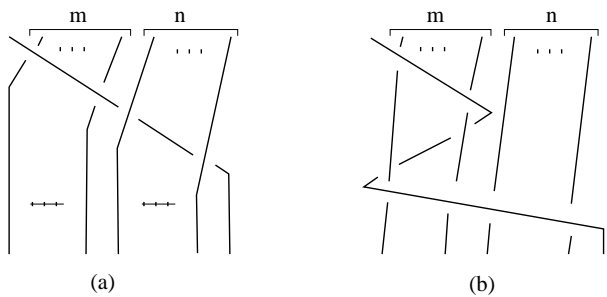}
\caption{The braids (a) $\beta_{m,n}$ and (b) $\sigma_{m,n}$.}
\label{SmnBmn-fig}
\end{center}
\end{figure}

It turns out that one may find smaller dilatations (or produce a reducible or periodic braid)
by passing a strand of $\beta_{m,n}$ 
once around the remaining strands.
As a particular example, we
consider the braids $\sigma_{m,n}$ defined by taking the rightmost-strand
of $\beta_{m,n}$ and passing it counter-clockwise once around the remaining strands.  
Figure~\ref{SmnBmn-fig} gives an illustration of $\beta_{m,n}$ and $\sigma_{m,n}$.
The braid $\sigma_{1,3}$ is conjugate to Ham and Song's braid given in (\ref{HamSong-braid})
(see Section~\ref{Graphmap-section}).

For $|m-n| \leq 1$, we show that $\sigma_{m,n}$ is periodic (case $m=n$) or
reducible (case $|m-n|=1$).    Otherwise $\sigma_{m,n}$ is pseudo-Anosov
with dilatation strictly less than the dilatation of $\beta_{m,n}$.
The dilatations of $\sigma_{g-1,g+1}$ satisfy the inequality 
\begin{eqnarray}\label{inequality-eqn}
\lambda(\sigma_{g-1,g+1})^g \leq 2 + \sqrt{3}.
\end{eqnarray}
(Section~\ref{Salem-Boyd-section}, Proposition~\ref{min-Smn-prop}).

Let $\sM_g^s$ denote the set of isotopy classes (called {\it mapping classes}) 
of a closed orientable genus $g$ surface $F$ fixing $s$ points.  For any subset
$\Gamma \subset \sM_g^s$, define $\lambda(\Gamma)$ to be the least dilatation
among pseudo-Anosov elements of $\Gamma$, and let $\delta(\Gamma)$ be
the logarithm of $\lambda(\Gamma)$.  For the braid group $\sB(D,s)$,
and any subset $\Gamma \subset \sB(D,s)$, define $\lambda(\Gamma)$ and
$\delta(\Gamma)$ in a similar way.
By a result of Penner \cite{Penner91} (see also, \cite{McMullen:Poly}), 
\begin{eqnarray}\label{Penner-asymp}
\delta(\sM_g^0) \asymp \frac{1}{g}. 
\end{eqnarray}

An element of $\sM_g = \sM_g^0$ is called {\it hyperelliptic} if it commutes with
an involution $\iota$ on $F$ such that the quotient of $F$ by $\iota$ is $S^2$. 
Let $\sM_{g,\mathrm{hyp}} \subset \sM_g$ denote the subset of {\it hyperelliptic} elements of
$\sM_g=\sM_g^0$.    Any pseudo-Anosov braid on $2g+1$ strands determines
a hyperelliptic element of $\sM_g$ with the same dilatation (see Proposition~\ref{spectrum-prop}).  Let
$\widehat{{\sB}}^{\mathrm{lift}}_{\mathrm{pA}}(D,2g+1) \subset \sM_{g,\mathrm{hyp}}$ denote the image
of these pseudo-Anosov braids in $\sM_{g,hyp}$.  Then
we have
$$
\delta(\sM_g) \leq \delta(\sM_{g,\mathrm{hyp}}) \leq \delta(\widehat{{\sB}}^{\mathrm{lift}}_{\mathrm{pA}}(D,2g+1)) = \delta(\sB(D,2g+1)).
$$
Putting this together with (\ref{inequality-eqn}) and (\ref{Penner-asymp}) we 
have the following.
\begin{thm}\label{asymp-thm}  For $g \ge 1$,
$$
\delta(\sB(D,2g+1)) \asymp \frac{1}{g}
$$
and
$$
\delta(\sM_{g,\mathrm{hyp}})  \asymp \frac{1}{g}.
$$
\end{thm}

The paper is organized as follows.  Section~\ref{prelim-section} reviews definitions
and properties of braids, mapping classes, and spectra of pseudo-Anosov dilatations.  In
Section~\ref{main-section}, we determine the Thurston-Nielsen classification of
$\beta_{m,n}$ and $\sigma_{m,n}$, by finding efficient and irreducible graph maps
for their monodromy actions following \cite{BH94}.  We observe that the
graph maps, and hence
the associated train tracks have ``star-like" components, and their essential forms don't
depend on $m$ and $n$ (see Figures~\ref{BmnGraph-fig} and \ref{SmnGraph-fig}).
To find bounds and inequalities among the dilatations, we apply the notion of 
Salem-Boyd sequences \cite{Boyd77}, \cite{Salem44}, and relate the similar forms
of the graph maps for $\beta_{m,n}$ and $\sigma_{m,n}$ to similar forms for
characteristic polynomials of the dilatations.
In particular, we show that
the least dilatation that occurs among $\beta_{m,n}$ and $\sigma_{m,n}$
for $m+n = 2g$ is realized by $\sigma_{g-1,g+1}$,
and find bounds for $\lambda(\sigma_{g-1,g+1})$ yielding the upper bound in (\ref{inequality-eqn}).

Section~\ref{discussion-section} discusses the problem of determining the least dilatations of 
special subclasses of  pseudo-Anosov maps. 
In Section~\ref{forcing-section}, we briefly describe the relation between the forcing relation 
on {\it braid types} and dilatations, and show how the $\sigma_{m,n}$ arise as the braid types of 
periodic orbits of the {\it Smale-horseshoe map}. 
In Section~\ref{fiberedQ-section}, we consider 
pseudo-Anosov maps arising as the monodromy of fibered links, and relate our examples to
those of Brinkmann.
\medskip

\noindent
{\sc Acknowledgements:}  The authors thank H. Minakawa for valuable 
discussions, and an
algebraic trick that improved our original upper bound for $\lambda(\sigma_{g-1,g+1})$.
The first author thanks the J.S.P.S., Osaka University and host Makoto Sakuma
for their hospitality and support during the writing of this paper.
The second author is grateful for  the financial support provided by the research fellowship of
the 21st century COE program in Kyoto University.

\section{Preliminaries}\label{prelim-section}

In this section, we review some basic definitions and properties of braids
(Section~\ref{braids-section}), mapping class groups (Section~\ref{mappings-section}),
 and spectra (Section~\ref{spectrum-section}).  The results are well-known, and more
 complete expositions can be found in \cite{Birman74}, \cite{CB88}, \cite{FLP}, and \cite{BH94}.
 We include them here for the convenience of the reader.

\subsection{Braids}\label{braids-section}

Let $F$ be a compact orientable surface with $s$ marked points 
${\sS}=\{p_1,\dots,p_s\} \subset \mathrm{int}(F)$.
A {\it braid representative} $\beta$ on $F$ is the union of $s$ {\it strands},  or images 
of continuous maps 
$$f_{p_1},\dots,f_{p_s} : I \rightarrow F \times [0,1],
$$
satisfying, for $i=1,\dots,s$,
\begin{description}
\item{(B1)} $f_{p_i} (0) = p_i$;
\item{(B2)} $f_{p_i}(1) \in \sS$; 
\item{(B3)} $f_{p_i}(t) \in F \times \{t\}$; and
\item{(B4)} $f_{p_i}(t) \neq f_{p_j}(t)$ for all $t$ and $i \neq j$.
\end{description}
Define the product of two braid representatives to be their
concatenation.  Let $\sB(F;\sS)$
be the set of braid representatives up to ambient isotopy fixing the 
boundary of $F$ point-wise.
The above definition of product determines a well-defined group structure on $\sB(F;\sS)$, 
and the group is called the {\it braid group} on $F$. 

For any partition $\sS = \sS_1 \cup \cdots \cup \sS_r$, let $\sB(F;\sS_1,\dots,\sS_r)$ be 
the subgroup of $\sB(F;\sS)$ consisting of braids $(f_{p_1},\dots,f_{p_s})$ satisfying:
for all $j=1,\dots,r$ and all $p \in \sS_j$,  
$$
f_p(1) \in \sS_j.
$$

\begin{figure}[htbp]
\begin{center}
\includegraphics[width=1.5in]{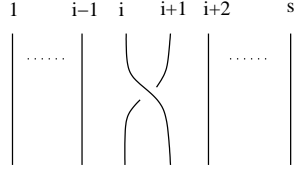}
\caption{The braid generator $\sigma_i$.}
\label{braidgen-fig}
\end{center}
\end{figure}

Let $F$ be either the closed disk $D$ or the 2-sphere $S^2$ and 
let $\sS \subset \mbox{int}(F)$ be a set of $s$ points.   The braid group $\sB(F;\sS)$ 
 has generators 
$\sigma_1,\dots,\sigma_{s-1}$, where $\sigma_i$ is the braid
 shown in Figure~\ref{braidgen-fig}.  
When $F = D$, the
group  $\sB(D;\sS)$ is called the {\it Artin braid group} and has finite presentation
$$
\langle \sigma_1,\dots,\sigma_{s-1}\ : 
\ \sigma_i\sigma_{i+1}\sigma_i = \sigma_{i+1}\sigma_i\sigma_{i+1}
\mbox{, and } \sigma_i\sigma_j = \sigma_j\sigma_i \ \mbox{if $|i-j| \ge 2$.}\rangle.
$$
\vspace{-12pt}
\vspace{-12pt}

\begin{figure}[htbp]
\begin{center}
\includegraphics[height=1.5in]{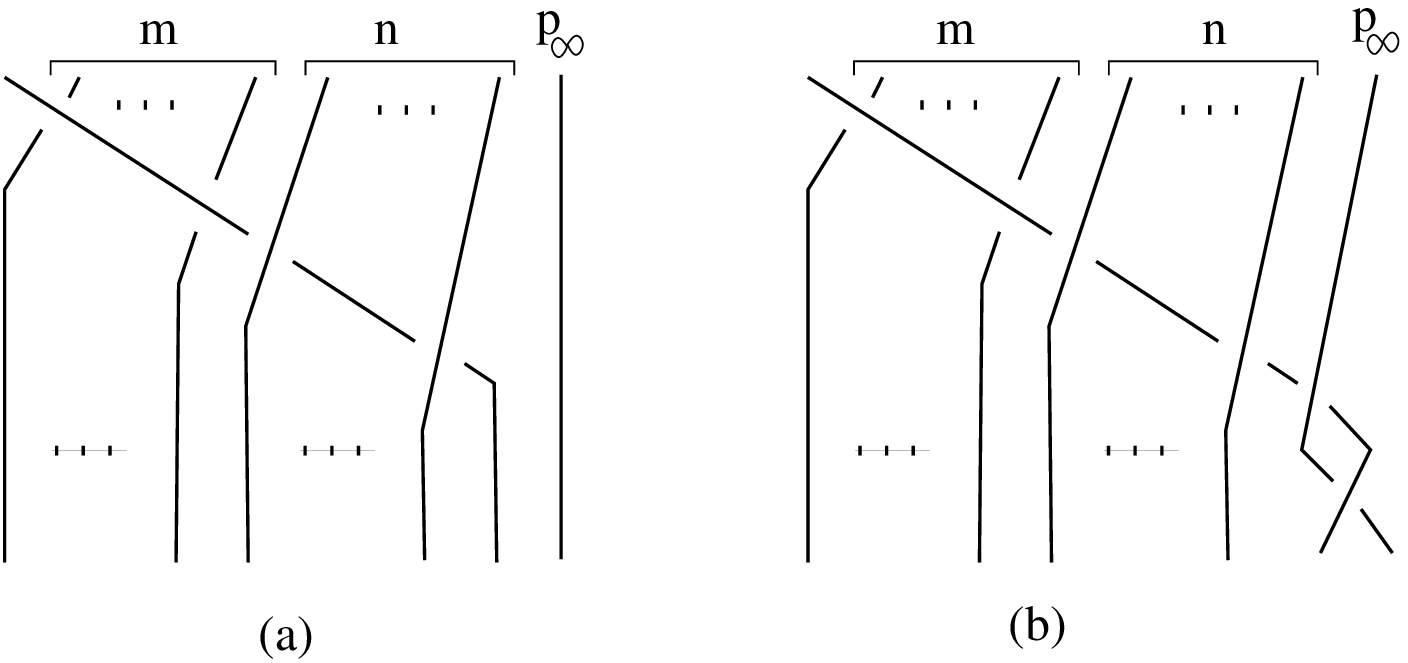}
\caption{Images of (a) $\beta_{m,n}$ and (b) $\sigma_{m,n}$ in $\sB(S^2;\sS,\{p_\infty\})$.}
\label{spherical-fig}
\end{center}
\end{figure}

Consider the natural map $c : D \rightarrow S^2$ given by identifying the boundary of $D$ to a
point $p_\infty$ on $S^2$.  By abuse of notation, we will write $\sS$ for $c(\sS)$.
Then there is an induced map:
\begin{eqnarray}\label{widehat-eqn}
\sB(D;{\sS}) &\rightarrow&\sB(S^2;\sS,\{p_\infty\})
\\
\beta &\mapsto& \widehat{\beta}.\nonumber
\end{eqnarray}
For example, $\widehat{\beta_{m,n}}$ and $\widehat{\sigma_{m,n}}$ are shown in 
Figure~\ref{spherical-fig} with the strand associated to $p_\infty$ drawn on the right.
For $\beta \in \sB(D;{\sS})$, let $\overline{\beta} \in \sB(S^2;{\sS})$ be the image of
$\widehat{\beta}$ under the forgetful map
\begin{eqnarray}\label{overline-eqn}
 \sB(S^2;\sS,\{p_\infty\}) &\rightarrow&  \sB(S^2;{\sS})\\
\widehat{\beta} &\mapsto& \overline{\beta}\nonumber
 \end{eqnarray}
obtained by   the strand associated to $p_\infty$.   
 The following lemma can be found in \cite{Birman74}.

\begin{lem} 
The map $\sB(D;\sS) \rightarrow \sB(S^2;\sS)$ 
given by composing the maps in (\ref{widehat-eqn}) and (\ref{overline-eqn}) has kernel 
normally generated by
$$
\xi=\xi_s = \sigma_1\sigma_2 \cdots \sigma_{s-1}^2\sigma_{s-2} \cdots \sigma_1.
$$
\end{lem}
For example, the braids $\beta_{m,n}$ and $\sigma_{m,n}$ shown in Figure~\ref{SmnBmn-fig}
differ by a conjugate of $\xi$ and hence  we have the following.

\begin{prop}  The braids $\beta_{m,n}$ and $\sigma_{m,n}$ satisfy
$$
\overline{\beta_{m,n}}= \overline{\sigma_{m,n}}.
$$
\end{prop}

The final lemma of this section deals with notation.  

\begin{lem}\label{invariance-lem} 
Let $F$ be $D$ or $S^2$. 
Let $\sS_1$ and $\sS_2$ be two finite subsets of $\mbox{int}(F)$ with the same cardinality, 
and $h : D \rightarrow D$ any homeomorphism taking $\sS_1$ to $\sS_2$.
Then conjugation by $h$ defines an isomorphism
$$
\sB(F;\sS_1) \rightarrow \sB(F;\sS_2).
$$
\end{lem}
\medskip

In light of Lemma~\ref{invariance-lem}
if $s$ is the cardinality of $\sS$, we will also write $\sB(F,s)$ for $\sB(F;\sS)$.

\subsection{Mapping class groups}\label{mappings-section}

For any closed orientable surface $F$ and a finite subset $\sS \subset F$ of {\it marked points}, 
let $\sM(F;\sS)$ be the group of isotopy classes (called {\it mapping classes}) of 
orientation preserving homeomorphisms of $F$ set-wise preserving 
$\sS$.

The Thurston-Nielsen classification states that any homeomorphism of a surface
is isotopic to one of three types, which we describe below.

A map $\Phi : F \rightarrow F$ is defined to be
{\it periodic} if
some power of $\Phi$ equals the identity map; and
 {\it reducible} if
there is a $\Phi$-invariant closed $1$-submanifold whose complementary components
in $F \setminus \sS$ have negative Euler characteristic.   A mapping class $\phi \in \sM(F;\sS)$
is {\it periodic} (respectively,  {\it reducible})
 if it contains a representative that is periodic (respectively, 
reducible).  

\begin{figure}[htbp]
\begin{center}
\includegraphics[height=1in]{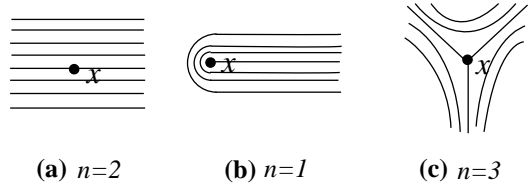}
\caption{Local picture of a singular foliation.}
\label{prongs-fig1}
\end{center}
\end{figure}

Before defining the  third type of mapping class, we will make some preliminary definitions.
A {\it singular foliation} $\sF$  on a closed surface $F$ with respect to the set of marked points $\sS$ 
is a partition of $F$ into a union of real intervals $(-\infty,\infty)$ and $[0,\infty)$ 
called {\it leaves}, such that,  for each point $x \in F$,  the foliation $\sF$ near $x$ has
(at least) one of the following types in a local chart around $x$:
\begin{description}
\item{(F1)} 
$x \in F$ is a {\it regular point} (we will also say a {\it $2$-pronged point}) of 
$\sF$ (Figure~\ref{prongs-fig1}(a));
\item{(F2)} 
$x \in F$ is an {\it $n$-pronged singularity of $\sF$} (Figure~\ref{prongs-fig1}(b) and \ref{prongs-fig1}(c)),
where
\begin{description}
\item{(i)} if $x \in  \sS$, $n \ge 1$,  and 
\item{(ii)} if $x \in F \setminus \sS$, $n \ge 3$.
\end{description}
\end{description}

Two singular foliations $\sF^{+}$ and $\sF^{-}$ with respect to $\sS$ are {\it transverse} if 
they have the same set of singularities 
and if  the leaves of $\sF^{+}$ and $\sF^{-}$ intersect transversally.

A path $\alpha$ on $F$ is a {\it transverse arc (relative to $\sF$)} if $\alpha$ intersects
the leaves of $\sF$ transversely.
 Two transverse arcs $\alpha_0$ and $\alpha_1$ are {\it relatively homotopic
  (with respect to $\sF$)} if
there is a homotopy $\alpha: [0,1] \times [0,1] \rightarrow F$ 
such that $\alpha ([0,1] \times \{0\}) = \alpha_0$, 
$\alpha ([0,1] \times \{1\}) = \alpha_1$ and, for all $a \in [0,1]$,  $\alpha(\{a\} \times [0,1])$ is contained in a leaf of ${\mathcal F}$.
We say that $\mu$ is a {\it transverse measure} 
on a singular foliation $\sF$ with respect to $\sS$  
if $\mu$ defines a non-negative Borel measure $\mu(\alpha)$ 
on each transverse arc $\alpha$  with the following two properties. 
\begin{description}
\item{(M1)} 
If $\alpha'$ is a subarc of $\alpha$, then $\mu(\alpha')=\mu(\alpha)|_{\alpha'}$.
\item{(M2)} 
If $\alpha_0$, $\alpha_1$ are  arcs transverse to ${\mathcal F}$, and relatively
homotopic with respect to ${\mathcal F}$, 
 then $\mu(\alpha_0) = \mu(\alpha_1) $. 
\end{description}

A pair $({\mathcal F}, \mu)$ satisfying (M1) and (M2) is called a {\it measured foliation}. 
Given a measured foliation $(\sF,\mu)$ on F with respect to $\sS$ and 
a number $\lambda>0$, $(\sF, \lambda \mu)$ denotes the measured foliation whose
leaves are the same as those of $\sF$ and such that, for all arcs $\alpha$,
the measure of $\alpha$  transverse to ${\mathcal F}$ 
is given by $\lambda \mu(\alpha)$. 
For a homeomorphism $f: F \rightarrow F$ set-wise preserving $\sS$, 
 $(\sF', \mu')= f({\mathcal F}, \mu)$ is the measured foliation whose
leaves  are the images 
of leaves of ${\mathcal F}$ under $f$, 
and the measure $\mu'$ on each arc $\alpha$ transverse to $\sF'$
is given by $\mu(f^{-1}(\alpha))$. 

The final classification type for surface homeomorphisms is defined as follows.
A map $\Phi : F \rightarrow F$ is  {\it pseudo-Anosov}, if 
there is a number $\lambda >1$ and a pair of transverse measured foliations 
$(\mathcal{F}^\pm, \mu_\pm)$ with respect to $\sS$ such that 
$\Phi(\mathcal{F}^\pm, \mu_\pm) = (\mathcal{F}^\pm, \lambda^{\pm 1} \mu_\pm)$.
The number $\lambda(\Phi) = \lambda$ is called the {\it dilatation} of $\Phi$, and 
$\mathcal{F}^+$ and $\sF^{-1}$ are called the {\it stable} and {\it unstable} foliations 
or the {\it invariant foliations} associated to $\Phi$.

A mapping class $\phi \in \sM(F;\sS)$ is {\it pseudo-Anosov} if $\phi$ is 
the isotopy class of a pseudo-Anosov map $\Phi$.  In this case, the dilatation of $\phi$ is
defined to be $\lambda(\phi) = \lambda(\Phi)$.

\begin{thm} [Thurston-Nielsen Classification Theorem]  Any element 
$$\phi \in \sM(F;\sS)$$
is either periodic, reducible or pseudo-Anosov.  
Furthermore, if $\phi$ is  pseudo-Anosov, then the pseudo-Anosov representative of
$\phi$ is unique.
\end{thm}

As with braids, for any partition $\sS = \sS_1 \cup  \cdots \cup \sS_r$, there is a subgroup
$$
\sM(F;\sS_1,\dots,\sS_r) \subset \sM(F;\sS)
$$
 that preserves each $\sS_i$ setwise.  
If $r > 1$, there is a natural map
$$
\sM(F;\sS_1,\dots,\sS_r) \rightarrow \sM(F;\sS_1,\dots,\sS_{r-1})
$$
called the {\it forgetful map}.  For pseudo-Anosov mapping classes $\phi$, 
$\log(\lambda(\phi))$  can be interpreted as the minimal topological 
entropy among all representatives of the class $\phi$ \cite{FLP}.  
We thus have the following inequality on dilatations.

 \begin{lem}\label{closure1-lem}  
 Let $\phi \in \sM(F;\sS_1,\dots,\sS_r)$ and let $\psi \in \sM(F;\sS_1,\dots,\sS_{r-1})$ be the image of $\phi$ 
 under the forgetful map.   If $\phi$ and $\psi$ are both pseudo-Anosov, then
 $$
 \lambda(\phi) \ge \lambda(\psi).
 $$
 \end{lem}

\begin{lem}\label{closure2-lem} 
Let $\phi \in \sM(F;\sS_1,\dots,\sS_r)$ be pseudo-Anosov. 
Suppose that the pseudo-Anosov representative $\Phi$ of $\phi$ does not have a 1-pronged
singularity at any point of $S_r$.  
Let $\psi \in \sM(F;\sS_1,\dots,\sS_{r-1})$ be the image of $\phi$ under the forgetful map. 
Then $\psi$ is pseudo-Anosov and $\lambda(\psi)= \lambda(\phi)= \lambda(\Phi)$. 
\end{lem}

\proof 
Let $\sF^{\pm}$ be singular foliations with respect to $\sS_1\cup\dots\cup\sS_r$, and let
$(\sF^\pm,\mu_\pm)$ be a pair of transverse measured foliations associated to $\Phi$. 
Since $\sF^{\pm}$ does not have 1-pronged singularities at points of $S_r$,  the
$\sF^{\pm}$ give well-defined singular foliations with respect to $\sS_1\cup\dots\cup\sS_{r-1}$.
Thus, $\Phi$ is a pseudo-Anosov representative of $\psi$ in $\sM(F;\sS)$, and 
hence $\lambda(\psi)= \lambda(\phi)= \lambda(\Phi)$. 
\qed
\medskip

As in the case of braids, changing the location
of the points in $\sS$ by a homeomorphism does not change the group
$\sM(F;\sS)$  (cf. Lemma~\ref{invariance-lem}).

\begin{lem}\label{invarianceM-lem}
Let $\sS_1$ and $\sS_2$ be two finite subsets of $F$ with the same cardinality
and let $h : F \rightarrow F$ be any homeomorphism taking $\sS_1$ to $\sS_2$.
Then conjugation by $h$ defines an isomorphism
$$
\sM(F;\sS_1) \rightarrow \sM(F;\sS_2).
$$
\end{lem}
\medskip

\noindent
If $F$ has genus
$g$, and $\sS$ has cardinality $s$, we will also write 
$\sM_g^s = \sM(F;\sS)$.

\begin{figure}[htbp]
\begin{center}
\includegraphics[width=2in]{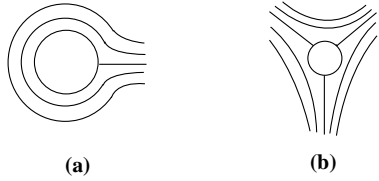}
\caption{Leaves of a singular foliation near a boundary component. }
\label{prongs-fig2}
\end{center}
\end{figure}

The theory of mapping class groups on closed surfaces extends to mapping class groups on
surfaces with boundary.  Let $F^b$ be a compact orientable surface with $b$ boundary 
components, and let $\sS \subset \mathrm{int}(F^b)$ be a finite set.  
Define $\sM(F^b;\sS)$ to  be the group of isotopy classes (or {\it mapping classes})
of orientation preserving homeomorphisms 
of $F^b$ set-wise preserving $\sS$ and the boundary components. 
A {\it singular foliation} $\sF$ on $F^b$ 
with respect to the set of marked points $\sS$ is a partition of $F$ into a union of leaves 
  such that each point $x \in \mbox{int}(F)$ has a local chart satisfying one of the conditions 
(F1), (F2) in the definition of the singular foliation on the closed surface, and 
each boundary component has $n$-prongs for some $n \ge 1$.  Figure~\ref{prongs-fig2} 
illustrates representative leaves of a singular foliation with a 1-pronged (Figure~\ref{prongs-fig2}(a))
and 3-pronged (Figure~\ref{prongs-fig2}(b)) singularity. 
{\it Periodic}, {\it reducible} and {\it pseudo-Anosov} mapping  classes in $\sM(F^b;\sS)$
are defined as for the case of closed surfaces using this definition of 
singular foliations. 

Let 
\begin{eqnarray}\label{contraction-map}
c : F^b \rightarrow \overline{F^b}
\end{eqnarray}
be the continuous map where $\overline{F^b}$ is the closed surface obtained
from $F^b$ by contracting the boundary components to points $q_1,\dots,q_b$.
As before, we will write $\sS$ for $c(\sS)$.
Let $\sQ = \{q_1,\dots,q_b\}$.   The above definitions imply the following. 

\begin{lem}\label{boundary-lem}  The contraction map $c$ in (\ref{contraction-map}),
induces an isomorphism
$$
c_* : \sM(F;\sS) \rightarrow \sM(\overline{F^b}; \sS,\sQ),
$$
which preserves the Thurston-Nielsen classification of mapping classes. 
Furthermore, if $\sF$ is a singular foliation defined on $F$ which is $n$-pronged along
a boundary component $A$ of $F$,  then the image of $\sF$ under
$c_*$ has an $n$-pronged singularity at $c_*(A)$. 
\end{lem}

\begin{figure}[htbp]
\begin{center}
\includegraphics[width=2in]{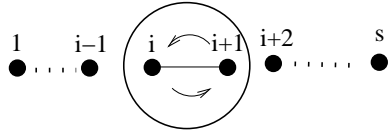}
\caption{Action of $\sigma_i$ as a homeomorphism of $F$.}
\label{braidmap-fig}
\end{center}
\end{figure}

\noindent
The isomorphism $c_*$ given in Lemma~\ref{boundary-lem} is handy in discussing
mapping  classes coming from braids.
Let $F$ be either the closed disk $D$ or the 2-sphere $S^2$.  
There is a natural homomorphism
\begin{eqnarray}\label{braidmonodromy-eqn}
\sB(F;\sS)& \rightarrow& \sM(F;\sS)\\
\beta &\mapsto& \phi_\beta \nonumber
\end{eqnarray}
defined as follows.  
Let $D_1,\dots,D_{s-1} \subset \mbox{int}(D)$ be disks with 
$D_i \cap D_j = \emptyset$ for $i \ne j$ 
such that $D_i$
contains two points $p_{i}$ and $p_{i+1}$ of $\sS$ and no other points of $\sS$.
The action of a generator $\sigma_i$ of $\sB(F;\sS)$
is the mapping class in $\sM(F;\sS)$ that fixes the exterior of $D_i$ 
and rotates a closed line segment connecting $p_{i}$ and $p_{i+1}$ in $D_i$ by 
$180$ degrees in the counter-clockwise direction as in Figure~\ref{braidmap-fig}.

Given a braid $\beta \in \sB(D;\sS)$, let
$\widehat{\beta}$ be its image in $\sB(S^2;\sS,\{p_\infty\})$ as in (\ref{widehat-eqn}).  
Then the isomorphism $c_*: \sM(D;\sS) \rightarrow \sM(S^2;\sS,\{p_\infty\})$ induced
by the contraction map
of Lemma~\ref{boundary-lem} satisfies
\begin{eqnarray}\label{cmap-eqn}
c_*(\phi_\beta) = \phi_{\widehat{\beta}}.
\end{eqnarray}

\noindent
The following useful lemma can be found in \cite{Birman74}.

\begin{lem}  
\label{braidcenter-lem}
If $\sS$ has cardinality $s$, then the kernel of the map
$$\sB(D;\sS) \rightarrow \sM(S^2;\sS,\{p_\infty\})$$ defined by
$\beta \mapsto \phi_{\widehat{\beta}}$ is the center $Z(\sB(D;\sS))$ 
generated by a full twist braid 
$$
\Delta=(\sigma_1\cdots \sigma_{s-1})^s.
$$
\end{lem}

We say that $\beta$ is {\it periodic} (respectively,  {\it reducible}, {\it pseudo-Anosov}),
 if $\phi_{\widehat{\beta}}$ is periodic (respectively,   reducible, pseudo-Anosov) as an element of 
 $\sM(S^2;\sS,\{p_\infty\})$.  
 In the pseudo-Anosov case, we set $\lambda(\beta) = \lambda(\widehat{\beta}) =  \lambda(\phi_{\widehat{\beta}})$.   
 
Let $\overline{\beta}$ be
the image of $\widehat{\beta}$ in $\sB(S^2;\sS)$ under the forgetful map  (\ref{overline-eqn}).  Then
Lemma~\ref{closure1-lem} implies that if $\widehat{\beta}$ and $\overline{\beta}$ are 
pseudo-Anosov, we have
\begin{eqnarray}\label{overline-ineq}
\lambda(\beta) = \lambda(\widehat{\beta}) \ge \lambda(\overline{\beta})
\end{eqnarray}
and by Lemma~\ref{closure2-lem}
equality holds if the point $p_{\infty}$ is not a $1$-pronged singularity for 
the invariant foliations associated to the pseudo-Anosov representative of $\phi_{\widehat{\beta}}$.

\subsection{The braid spectrum}\label{spectrum-section}

For any subset $\Gamma \in \sM_g^s$, let $\Sigma(\Gamma)$ be the set of logarithms
of dilatations coming from pseudo-Anosov elements of $\Gamma$.   Let
$\widehat{\sB}(D,s)$ be the image of $\sB(D,s)$ in $\sM_0^{s+1}$, 
and ${\widehat{\sB}}_{\mathrm{pA}}(D,s)$ the pseudo-Anosov elements of
$\widehat{\sB}(D,s)$.

\begin{prop}\label{spectrum-prop} For $g \ge 1$,
$$
\Sigma (\widehat{\sB}(D,2g+1)) \subset \Sigma(\sM_{g,\mathrm{hyp}} \subset \Sigma(\sM_g^0).
$$
\end{prop}

\proof  
Let $\sS \subset D$ be a subset of $2g+1$ points, and let $\widehat{\sS} = \sS \cup \{p_\infty\}$.
Let $F$ be the double cover of $S^2$ branched along  $\widehat{\sS}$.   Then $F$ has
genus $g$.  We will define a set map
$$
\widehat{{\sB}}_{\mathrm{pA}}(D,2g+1) \rightarrow \sM(F,0)
$$
whose image consists of hyperelliptic elements, and which preserves dilatation.

Let $\phi \in \widehat{{\sB}}_{\mathrm{pA}}(D,2g+1)$.  Then $\phi$ has a unique
pseudo-Anosov representative homeomorphism $\Phi$.   Let $\Phi'$ be its lift to $F$
by the covering $F \rightarrow S^2$, with stable and unstable invariant foliations 
given by the lifts of the invariant foliations associated to $\Phi$.   Then $\Phi'$ is
pseudo-Anosov with the same dilatation as $\Phi$.  Let $\phi'$ be its isotopy class.
Then $\phi'$ defines a hyperelliptic, pseudo-Anosov
mapping class in $\sM(F;{\widehat{\sS}}')$, where 
${\widehat{\sS}}'$ is the preimage of $\widehat{\sS}$ in $F$, with the same dilatation
as $\phi$.

Now consider the forgetful map $\sM(F_g;\sS') \rightarrow \sM(F_g;\emptyset)$.  
The stable and unstable foliations associated to $\Phi'$ 
have prong orders at $\widehat{\sS}'$ that are divisible by the degree of the covering
$F \rightarrow S^2$  Thus, the singularities of $\Phi'$ at $\widehat{\sS}'$ are all 
even-pronged.
It follows that by Lemma~\ref{closure2-lem}, the image of $\phi'$ under the forgetful map
is pseudo-Anosov and has the same dilatations as $\phi'$.
\qed
\medskip
\\
Proposition~\ref{spectrum-prop} immediately implies the following corollary.

\begin{cor}\label{spectrum-cor}
$$
\delta (\sM_g^0) \leq \delta (\sM_{g,\mathrm{hyp}}) \leq  \delta (\sB(D,2g+1)).
$$
\end{cor}

\subsection{Criterion for the pseudo-Anosov property}\label{Bestvina-Handel-section}

What follows is a criterion for determining when a braid $\beta \in \sB(D,s)$ is pseudo-Anosov 
\cite{BH94} (see also \cite{Penner:pA} and \cite{PH:TrainTracks}). 

Let $G$ be a finite graph embedded on an orientable surface $F$, possibly with self-loops, but 
no vertices of valence $1$ or $2$.  Let 
 $\DE(G)$ be the  set of directed edges of $G$.
 For $e \in \DE(G)$,  let $s(e)$ and $t(e)$ be the
initial ({\it source}) vertex and end ({\it target}) vertex respectively,
and let $\overline{e}$ be the same edge with opposite orientation.
An {\it edge path} $\tau$ is a directed path $\tau=e_1\cdots e_\ell$
on $G$ where $e_1,\dots,e_\ell \in \DE(G)$
satisfies $t(e_i) = s(e_{i+1})$ for $i=1,\dots,\ell-1$.  

\begin{figure}[htbp]
\begin{center}
\includegraphics[width=3in]{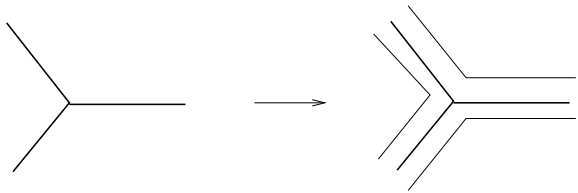}
\caption{Fat graph near a trivalent vertex.}
\label{fatgraph-fig}
\end{center}
\end{figure}

Let   $\sF(G)$ be the {\it fat graph} of $G$ obtained by
taking a regular neighborhood of the embedded graph $G$ on $F$,
(see Figure~\ref{fatgraph-fig}).    Let $\proj : \sF(G) \rightarrow G$ be
the natural projection to $G$.
An {\it (embedded) edge path} $\tau$ is a path on $\sF(G)$ on which $\proj$
is locally injective, and whose endpoints lie on vertices of $G$.  
Let $G$ and $H$ be two embedded graphs on the surface $F$.  
A  {\it graph map} 
$$
\g : G \rightarrow \sF(H)
$$
is a map sending each edge path on $G$ to an embedded edge path on $\sF(H)$,
so that if an edge path $\tau$ equals the product of edge paths $\tau_1$ and
$\tau_2$ on $G$, then $\g(\tau) = \g(\tau_1)\g(\tau_2)$ is a product of embedded
edge paths on $H$. 
Define the {\it composition} $\g\h : G \rightarrow \sF(H)$
 of two graph maps
  $\h : G \rightarrow \sF(H)$ and $\g : H \rightarrow \sF(K)$, by
 $\g(\h(e)) = \g(\proj(\h(e))$. 
 When $K=H$, write $\g^k$ for the composition of $\g$ with itself $k$ times.

Consider the set of undirected edges $\sE^{\mathrm{tot}}(G)$ of $G$ and let $V^{\mathrm{tot}}(G)$ be the vector space 
of formal sums
$$
\sum_{i=1}^n a_i e_i,
$$
where $a_i \in \R$ and $e_i \in \sE^{\mathrm{tot}}(G)$.  Any edge path on $G$ determines an element of $V(G)$  by 
treating each directed edge as an undirected edge with coefficient 1, regardless of
orientation.
For a graph map $\g : G \rightarrow \sF(H)$, 
define the associated {\it  total transition matrix}  for $\g$ to be the transformation
$$\sT_\g^{\mathrm{tot}} : V^{\mathrm{tot}}(G) \rightarrow V^{\mathrm{tot}}(H)$$
taking each $e \in \sE^{\mathrm{tot}}(G)$ to 
$\proj(\g(e))$ considered as an element of $V^{\mathrm{tot}}(H)$.  

For the rest of this section, we restrict to the case $F = D$.
Let 
$$
{\sS}= \{p_1,\dots,p_s\} \subset \mbox{int}(D)$$
 be a set of marked points,
let $P_i$ be a small circle centered at $p_i$ whose interior disk does not contain any
other points of $\sS$, and let 
$$
P= \displaystyle\bigcup_{i=1}^s P_i.
$$
Choose a finite graph $G$ embedded on $D$ 
that is homotopy equivalent to $D \setminus \sS$ such that 
each circle in $P$ is a subgraph of  $G$.  
Given $\beta \in \sB(D,s)$, the mapping class $\phi_{\beta} \in \sM(D; \sS)$ 
induces a graph map $\g: G \rightarrow \sF(G)$ set-wise preserving $P$ 
(i.e, $\pi \g(P)=P$) and a total transition matrix $\sT_\g^{\mathrm{tot}}$.
Let pre$P$ be the set of edges $e$ of $G$ such that $\g^k(e)$ is contained in $P$ for
some $k \ge 1$.

By the definition of $P$ and pre$P$, ${\sT}_\g^\mathrm{tot}$ has the following form: 
$$
{\sT}_\g^{\mathrm{tot}}= 
\left(\begin{array}{ccc}\mathcal{ P}& \mathcal{A} & \mathcal{B} \\0 & \mathcal{Z} & \mathcal{C} \\0 & 0 & \sT\end{array}\right), 
$$
where $\mathcal{P}$ (respectively,  $\mathcal{Z}$) are the transition matrices  associated to 
$P$ (respectively,   pre$P$), and 
$\sT$ is the transition matrix associated to the rest of edges $\sE^{\mathrm{real}}(G)$.
Let $V(G)$ be the subspace of $V^{\mathrm{tot}}(G)$ spanned
by $\sE(G) = \sE^{\mathrm{real}}(G)$.
The matrix $\sT$ is the restriction of ${\sT}_\g^\mathrm{tot}$ to $V(G)$ and is called the {\it transition matrix
 (with respect to the real edges)}
 of $\phi_{\beta}$, associated to $\g$.  The spectral radius of $\sT$ is denoted by $\lambda(\sT)$.

Given a graph map $\g : G \rightarrow \sF(G)$, define the {\it derivative} of $\g$ 
$$
D_\g : \DE(G) \rightarrow \DE(G)
$$
as follows: for  $e \in \DE(G)$, write
$$
\proj(\g(e)) = e_1e_2 \cdots e_\ell, 
$$
where $e_i \in \DE(G)$, and set
$D_\g(e)$ to be the initial edge $e_1$.
A graph map $\g : G \rightarrow \sF(G)$ is 
{\it efficient} if for any $e \in \DE(G)$ and any $k \ge  0$,
$$
\proj(\g^k(e)) = e_1\cdots e_j
$$
where $D_\g(\overline{e_i}) \neq D_\g(e_{i+1})$, for all $i=1,\dots,j-1$.
We also say in this case that $\g^k$ has no {\it back track} for any $k \ge 0$.

A square matrix  $M$ is {\it irreducible}
if for every set of indices $i,j$, there is an integer $n_{i,j} > 0 $ such that the $(i,j)$-th entry of 
$M^{n_{i,j}}$ is strictly positive.  The following criterion for the pseudo-Anosov property
can be found in \cite{BH94}.

\begin{thm}
\label{BH-thm}
Let $\beta \in \sB(D,s)$, and let $\g:G \rightarrow \sF(G)$ be an induced graph map for $\phi_{\beta}$.   
Suppose that 
\begin{description}
\item{(BH:1)} $\g$ is efficient, and 
\item{(BH:2)} the associated transition matrix $\sT$ (with respect to the real edges) is 
irreducible with spectral radius greater than one.
\end{description}
Then $\beta$ is pseudo-Anosov with dilatation equal to $\lambda(\sT)$.
\end{thm}

It is not hard to check that the criterion of Theorem~\ref{BH-thm} behaves well under conjugation
of maps.  For the case of braids, this yields the following.

\begin{lem}\label{conj-lem}  Suppose $\alpha_1$ and $\alpha_2$ are conjugate braids with
$\alpha_2 = \gamma \alpha_1 \gamma^{-1}$. 
Then 
$\alpha_1$ is pseudo-Anosov if and only if $\alpha_2$ is pseudo-Anosov.  In 
this case, we have $\lambda(\alpha_1) = \lambda(\alpha_2)$ and a graph map
satisfying (BH:1) and (BH:2)  for $\phi_{{\alpha_2}}$ can be obtained from
one for $\phi_{{\alpha_1}}$ by the homeomorphism of $D$ corresponding to $\gamma$.
\end{lem}

\begin{figure}[htbp]
\begin{center}
\includegraphics[angle=0,height=1.3in]{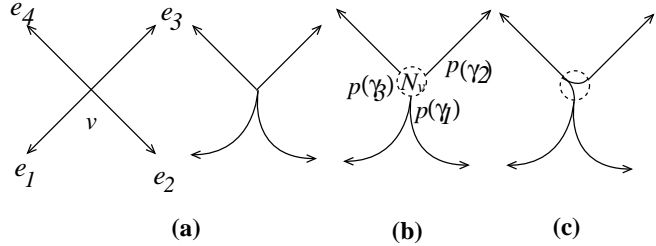}
\caption{Example of a graph smoothing.}
\label{fig_gate}
\end{center}
\end{figure}

Let $\phi \in \sM(D;\sS)$ be a pseudo-Anosov element, and let $\g : G \rightarrow \sF(G)$
be a graph map satisfying satisfying (BH:1) and (BH:2).  We construct an
{\it associated train track}
obtained by {\it graph smoothing} given   as follows.
Let $\mathcal{E}_v \subset \DE(G)$ be the set of oriented edges of $G$
emanating from a vertex $v$.  For $e_1,e_2 \in \mathcal{E}_v$,
$e_1$ and $e_2$  are  {\it equivalent}  if $D_\g^k(e_1)= D_\g^k(e_2)$ for some $k \ge 1$. 
A {\it gate} is an  equivalence class in $\sE_v$. 
The  train track $\tau_\g$ associated to $\g$ is constructed using the following steps: 
\begin{description}
\item{Step T1:}
Deform each pair of equivalent edges $e_1,e_2 \in \mathcal{E}_v$ in a small neighborhood of $v$ 
so that $e_1$ and $e_2$ are tangent at $v$, see Figure~\ref{fig_gate}(a). 
\item{Step T2:}
Insert  a small disk $N_v$ at each vertex $v$. 
For each gate $\gamma$, assign a point $p(\gamma)$ on the boundary of $N_v$, 
see Figure~\ref{fig_gate}(b). 
\item{Step T3:} If, for some edge $e$ of $G$ and some $k \ge 1$, 
$\g^k(e)$ contains consecutive edges $\overline{e_1'}e_2'$, 
with $e_1',e_2' \in {\mathcal E}_v$ and
such that $\gamma_1= [e_1']$ and $\gamma_2= [e_2']$, then 
join $p(\gamma_1)$ and $p(\gamma_2)$ by a smooth arc in $N_v$ 
which intersects the boundary of $N_v$ transversally at $p(\gamma_1)$
and $p(\gamma_2)$, and such that no two such arcs intersect in the
interior of $N_v$.  
\end{description}

For example, let $v$ be the initial point of 4 edges $e_1,\dots,e_4$.  
Assume there are 
three gates $\gamma_1=[e_1]=[e_2]$, $\gamma_2=[e_3]$ and $\gamma_3 = [e_4]$,
and that there are edges $f_1$ and $f_2$ of $G$ such that 
\begin{eqnarray*}
\g^r(f_1) &=& \cdots \overline{e_2}e_4 \cdots,\\  
\g^s(f_2) &=& \cdots \overline{e_3}e_4\cdots,
\end{eqnarray*}
for some $r,s \ge 1$.
Then Figure~\ref{fig_gate}(a) shows Step T1
applied to the edges $e_1$ and $e_2$;   Figure~\ref{fig_gate}(b) shows Step T2
applied to the gates $\gamma_1, \gamma_2$ and $\gamma_3$; and Figure~\ref{fig_gate}(c)
shows Step T3, which yields arcs connecting $p(\gamma_1)$ to $p(\gamma_3)$, and
$p(\gamma_2)$ to $p(\gamma_3)$.  

Let $v$ be a vertex of $G$.
The arcs constructed in Step T3 
are called {\it infinitesimal edges},
and the points $p(\gamma)$, which join two infinitesimal edges are 
called {\it cusps} of the train track.

If $\phi \in \sM(D;\sS)$ is pseudo-Anosov, and $\g: G \rightarrow G$ is an associated 
graph map satisfying (BH:1) and (BH:2), then the associated train track $\tau_\g$ 
constructed above
determines the invariant foliations $\sF^\pm$ 
associated to the pseudo-Anosov representative $\Phi$
of $\phi$.
(An explicit construction of $\Phi$ from 
the train track can be found in \cite{BH94}.)
In particular,  the number of prongs at the singularities of $\sF^\pm$ and at the
boundary of $D$ can be found in terms of the train track.

Each connected component $A$ of $D \setminus \tau_\g$ is either homeomorphic to an open
disk, or is a half-open annulus, whose boundary is the boundary of $D$.  In the former
case, the boundary of the closure of the connected component is a finite union of edges 
and vertices of $\tau_\g$.   If two of these edges meet at a cusp, then that cusp is said to 
 {\it belong to} $A$.    In the latter case, the closure of $A$ has two boundary components.
The boundary component which is not the boundary of $D$
is a finite union of edges and vertices of $\tau_\g$, and if two of these edges meet at
a cusp, we call the cusp an {\it exterior cusp} of $\tau_\g$.

\begin{lem}\label{prongs-lem} There is one singularity of $\sF^\pm$ in each of the open
complementary components of $\tau_\g$.  If $A$ is one such component, then the
singularity of $\sF^\pm$ in $A$ is $k$-pronged where $k$ is the number of cusps
of $\tau_\g$ belonging to $A$.  If $A$ contains the boundary of $D$, then the 
boundary of $D$ is $k$-pronged,
where $k$ is the number of exterior cusps of $\tau_\g$.
\end{lem}

\section{Main examples}\label{main-section}

This section contains properties of the braids $\beta_{m,n}$ and $\sigma_{m,n}$.  In Section~\ref{symmetry-section}, we show that the Thurston-Nielsen classifications of
$\beta_{m,n}$ and $\sigma_{m,n}$ do not depend on the order of $m$ and $n$.
In Section~\ref{Graphmap-section}, we find the Thurston-Nielsen classification
of the braids $\beta_{m,n}$ and $\sigma_{m,n}$, and in Section~\ref{dilatations-section},
we compute their dilatations in the pseudo-Anosov cases.
Section~\ref{TT-section} gives the train tracks for $\phi_{\beta_{m,n}}$ and $\phi_{\sigma_{m,n}}$
In Section~\ref{Salem-Boyd-section}, 
we use Salem-Boyd sequences to study inequalities between and the asymptotic behavior
of the dilatations of $\beta_{m,n}$ and $\sigma_{m,n}$.

\subsection{Symmetries of $\beta_{m,n}$ and $\sigma_{m,n}$}\label{symmetry-section}

\begin{figure}[htbp]
\begin{center}
\includegraphics[width=5in]{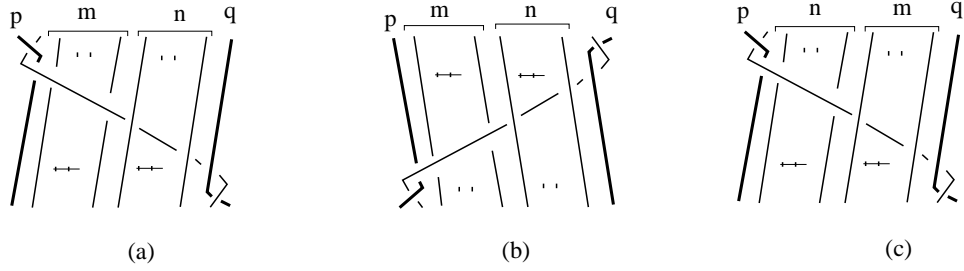}
\caption{Symmetry of $\beta_{m,n}^+$. }
\label{Bmn-sym-fig}
\end{center}
\end{figure}

Consider the  braid $\beta_{m,n}^+ \in \sB(D;\sS,\{p\},\{q\})$ drawn in Figure~\ref{Bmn-sym-fig}(a).

\begin{lem}\label{Bmnplus-conj-lem} The braid $\beta_{n,m}^+$ is conjugate to the
inverse of $\beta_{m,n}^+$.
\end{lem}

\proof
The inverse of $\beta_{m,n}^+$ is drawn in Figure~\ref{Bmn-sym-fig}(b).
Assume without loss of generality that the points $\sS \cup \{p\}\cup\{q\}$ are 
evenly spaced along a line $\ell$.
Let $\eta \in \sB(D;\sS\cup\{p,q\})$ be the braid obtained by
a half-twist of $\ell$ around the baricenter of $\sS \cup \{p\}\cup\{q\}$.  Then
conjugating the inverse of $\beta_{m,n}^+$ by $\eta$ in the larger group
$\sB(D;\sS\cup\{p,q\})$ yields
$\beta_{n,m}^+$ shown in Figure~\ref{Bmn-sym-fig}(c).
\qed
\medskip

\begin{lem}\label{Bmnplus-sym-lem} The braid $\beta_{m,n}$ is the image of
 $\beta_{m,n}^+$
under the forgetful map
$$
\sB(D;\sS,\{p\},\{q\}) \rightarrow \sB(D;\sS),
$$
and hence $\beta_{n,m}$ is conjugate to $\beta_{m,n}^{-1}$.
\end{lem}

\proof Compare Figure~\ref{Bmn-sym-fig}(a) with Figure~\ref{spherical-fig}(a) to get the first
part of the claim.  Since homomorphisms preserve inverses and conjugates, the
rest follows from Lemma~\ref{Bmnplus-conj-lem}. \qed
\medskip

\begin{lem}\label{Bmnplus-lem}  The mapping class $\phi_{\beta_{n,m}}$  is conjugate to $\phi_{\beta_{m,n}}^{-1}$.
\end{lem}

\proof The mapping class $\phi_{\beta_{m,n}}$ is the image of $\phi_{\beta_{m,n}^+}$
under the forgetful map
$$
\sM(D;\sS,\{p\},\{q\}) \rightarrow \sM(D;\sS).
$$
The rest  follows from 
Lemma~\ref{Bmnplus-sym-lem}.\qed
\medskip

\begin{prop}\label{Bmn-conj-prop} The Thurston-Nielsen classification of
 $\beta_{n,m}$ is the same as that of
$\beta_{m,n}$.
\end{prop}

\proof The Thurston-Nielsen classification of a mapping class is preserved under
inverses and conjugates.  Thus, the claim follows
from Lemma~\ref{boundary-lem} and  Lemma~\ref{Bmnplus-lem}. \qed
\medskip

\begin{figure}[htbp]
\begin{center}
\includegraphics[height=1.25in]{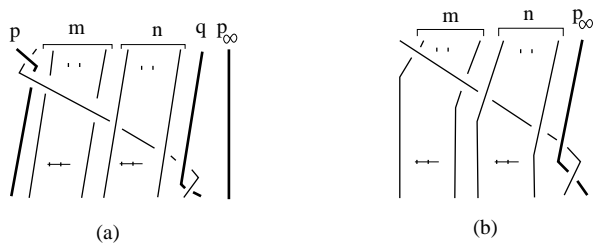}
\caption{Spherical braids (a) $\widehat{\beta_{m,n}^+}$ and (b) $\widehat{\sigma_{m,n}}$. }
\label{SphericalBmnplus-sym-fig}
\end{center}
\end{figure}

\begin{figure}[htbp]
\begin{center}
\includegraphics[height=0.75in]{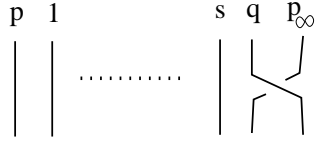}
\caption{Switching the roles of $p_\infty$ and $q$. }
\label{ConjBraid-fig}
\end{center}
\end{figure}

We now turn to $\sigma_{m,n}$. 
Let $\widehat{\beta_{m,n}^+}$ be the spherical braid associated to $\beta_{m,n}^+$
drawn in Figure~\ref{SphericalBmnplus-sym-fig}(a) and let 
$\nu \in \sB(S^2;\sS,\{p\},\{q,p_\infty\})$
be the spherical braid drawn in Figure~\ref{ConjBraid-fig}.
\vspace{-12pt}

\begin{lem}\label{BmnSmn-lem} The spherical braid $\widehat{\sigma_{m,n}}$ is the image of 
$\nu\widehat{\beta_{m,n}^+}\nu^{-1}$ under the forgetful map
$$
\sB(S^2;\sS,\{p\},\{q\},\{p_\infty\}) \rightarrow \sB(S^2;\sS,\{p_\infty\}), 
$$
and hence $\widehat{\sigma_{n,m}}$ is
conjugate to $\widehat{\sigma_{m,n}}^{-1}$.
\end{lem}

\proof Compare Figures~\ref{SphericalBmnplus-sym-fig}(a) and \ref{SphericalBmnplus-sym-fig}(b).  \qed

\begin{rem} In the statement of Lemma~\ref{BmnSmn-lem}, $\nu$ could be replaced 
by any braid which is the identity 
on $p$ and $\sS$, and interchanges $q$ and $p_\infty$.
\end{rem}

\vspace{-12pt}

\begin{cor}\label{BmnSmn-cor} The mapping class $\phi_{\sigma_{n,m}}$ is conjugate
to $\phi_{\sigma_{m,n}}^{-1}$.
\end{cor}

\proof
The mapping class $\phi_{\widehat{\sigma_{m,n}}}$ is the image of 
$\phi_\nu \phi_{\widehat{\beta_{m,n}^+}} \phi_\nu^{-1}$ under the forgetful map
$$
\sM(S^2;\sS,\{p\},\{q\},\{p_\infty\})\rightarrow \sM(S^2;\sS,\{p_\infty\}),
$$
and hence $\phi_{\widehat{\sigma_{n,m}}}$ is conjugate to $\phi_{\widehat{\sigma_{m,n}}}^{-1}$.
Since the contraction map $c$ defined in Equation (\ref{contraction-map}) induces an isomorphism
on mapping class groups the claim follows.
\qed

\begin{prop}\label{Smn-sym-prop}  The Thurston-Nielsen classification of $\sigma_{m,n}$
is the same as that for $\sigma_{n,m}$.
\end{prop}

\proof By Lemma~\ref{boundary-lem}, the Thurston-Nielsen classification for
$\phi_{\sigma_{m,n}}$ is the
same as that for $\phi_{\widehat{\sigma_{m,n}}}$.
The rest of the proof is similar to that of Proposition~\ref{Bmn-conj-prop}.\qed
\medskip

\subsection{Graph maps}\label{Graphmap-section}

In this section, we determine the Thurston-Nielsen classifications
of the braids $\beta_{m,n}$ and $\sigma_{m,n}$.   

\begin{thm} \label{Bmn-classification-thm}  The braid $\beta_{m,n}$ is
pseudo-Anosov for all $m,n \ge 1$, and 
$$
\lambda(\beta_{m,n}) = \lambda(\beta_{n,m}).
$$
\end{thm}

\begin{figure}[htbp]
\begin{center}
\includegraphics[width=5in]{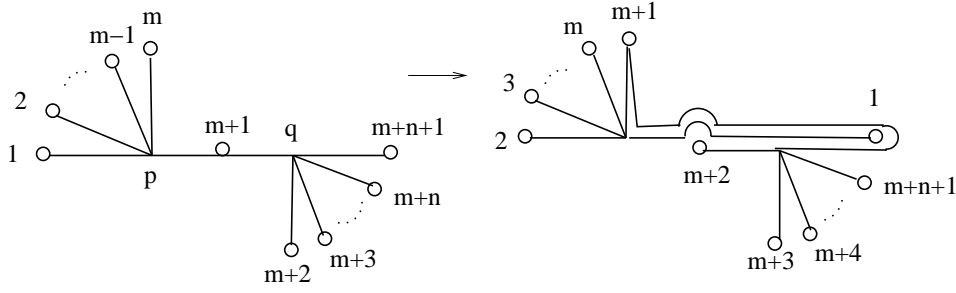}
\caption{Graph map $\g_{m,n}$ for $\phi_{\beta_{m,n}}$.}
\label{BmnGraph-fig}
\end{center}
\end{figure}

Consider the graph map $\g= \g_{m,n} : G_{m,n} \rightarrow \sF(G_{m,n})$  given in 
Figure~\ref{BmnGraph-fig}, where the ordering of the loop edges of the
graph $G_{m,n}$ corresponds to the left-to-right ordering of $\beta_{m,n}$.
Given a graph and a pair of vertices $a$ and $b$, let $e(a,b)$ denote the edge 
from vertex $a$ to vertex $b$.  
\vspace{-12pt}

\begin{prop}\label{BmnGraph-prop}
 The graph map $\g_{m,n} : G_{m,n} \rightarrow \sF(G_{m,n})$ is a graph map for
$\phi_{\beta_{m,n}}$ satisfying (BH:1) and (BH:2).
\end{prop}

\proof As shown in Figure~\ref{BmnGraph-fig}, any back track must occur at 
$e(p,m)$, that is, if $\g^k$ back tracks, and $k$ is chosen minimally,
then there is an edge $e \in \DE(G_{m,n})$ such that
\begin{eqnarray}\label{badedge-eqn}
\g^{k-1}(e) = \cdots \overline{e_1}\cdot e_2 \cdots
\end{eqnarray}
where
$$
D_\g(e_1) = D_\g(e_2) = e(m,p).
$$ 
This implies that $\overline{e_1} = e(p,m+1)$ and $e_2 = e(m+1,q)$ 
(or $\overline{e_1}=e(q,m+1)$ and $e_2=e(m+1,p)$).    As can be seen by
Figure~\ref{BmnGraph-fig}, one can verify that
 there can be no edge of the form given in (\ref{badedge-eqn}).
 This proves (BH:1).  
 
To prove (BH:2), it suffices to note that $\g^{m+n}(e(q,m+1))$ crosses all non-loop edges
 of $G_{m,n}$ in either direction, and for any non-loop edge $e$ of $G_{m,n}$, $\g^k(e)$ crosses $e(q,m+1)$ in either direction for some $k \leq \max\{m+1,n\}$. 
 \qed
\medskip

\noindent
{\bf Proof of Theorem~\ref{Bmn-classification-thm}.} By Proposition~\ref{BmnGraph-prop}, the braids $\beta_{m,n}$ are
pseudo-Anosov for all $m,n \ge 1$.  The claim now follows from
Proposition~\ref{Bmn-conj-prop}.\qed
 \medskip

We now turn to $\sigma_{m,n}$.
 \begin{figure}[htbp]
\begin{center}
\includegraphics[height=1.25in]{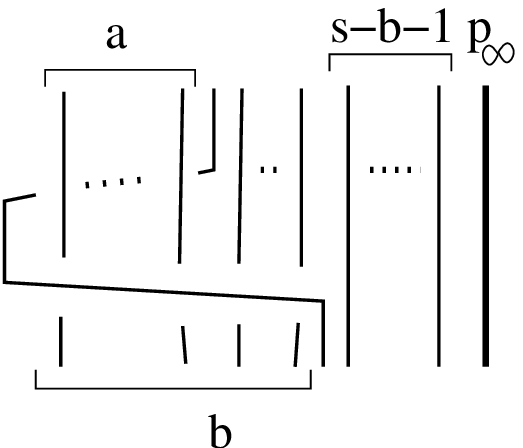}
\caption{Braid $\tau_{a,b} \in \sB(S^2;\sS,\{p_\infty\})$.}
\label{wrapbraid-fig}
\end{center}
\end{figure}

\begin{thm} \label{Smn-classification-thm}  The braid $\sigma_{m,n}$ is 
pseudo-Anosov for all $m,n \ge 1$ satisfying $|m-n| \ge 2$.  In these cases
$$
\lambda(\sigma_{m,n}) = \lambda(\sigma_{n,m}).
$$
For any $m \ge 1$,  $\sigma_{m,m}$ is
periodic, and $\sigma_{m,m+1}$ and $\sigma_{m+1,m}$ are
reducible.
\end{thm}

 In light of Proposition~\ref{Smn-sym-prop}, we will consider only $\sigma_{m,n}$ where
 $n \ge m \ge 1$.   To prove Theorem~\ref{Smn-classification-thm}, we 
 first redraw the braid $\sigma_{m,n}$ in a conjugate form
 using induction. 
 Let $\tau_{a,b}$ be the spherical braid drawn in Figure~\ref{wrapbraid-fig}.  
  Roughly speaking, conjugation by
 the braid $\tau_{a,b}$ on  $\sigma_{m,n}$
  is the same as passing a strand counterclockwise around the other
 strands, and then compensating below after a shift of indices.

 \begin{figure}[htbp]
\begin{center}
\includegraphics[width=4in]{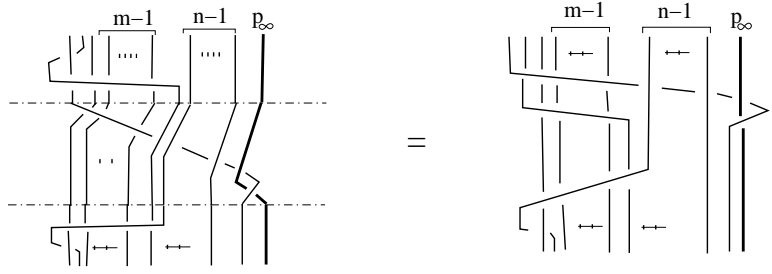}
\caption{Conjugating $\widehat{\sigma_{m,n}}$: initial step.}
\label{firstwrap-fig}
\end{center}
\end{figure}
\begin{figure}[htbp]
\begin{center}
\includegraphics[width=4.75in]{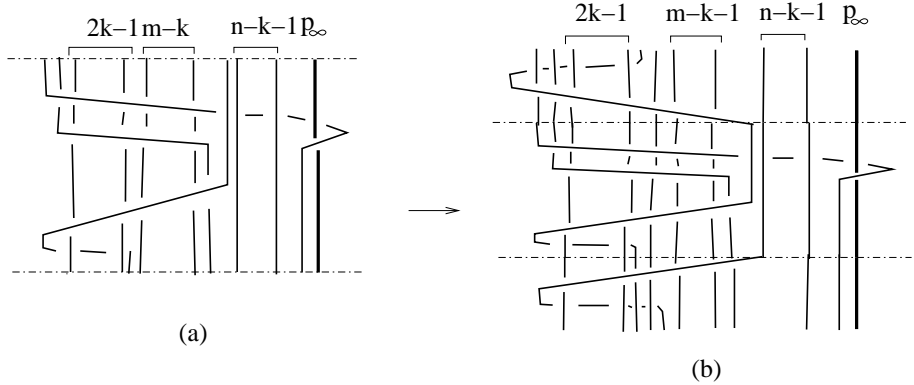}
\caption{Induction step.}
\label{inductive-fig}
\end{center}
\end{figure}

 Let ${\widehat{\sigma_{m,n}}}^{(0)} = \widehat{\sigma_{m,n}}$ be the image of $\sigma_{m,n}$ in $\sB(S^2;\sS,\{p_\infty\})$
as drawn in Figure~\ref{SphericalBmnplus-sym-fig}(b). 
Let 
$$
{\widehat{\sigma_{m,n}}}^{(1)} = \tau_{1,m+1} \widehat{\sigma_{m,n}} \tau_{1,m+1}^{-1}, 
$$
shown in Figure~\ref{firstwrap-fig}.  The 
inductive step is illustrated in Figure~\ref{inductive-fig}.  The $k$th braid $\widehat{\sigma_{m,n}}^{(k)}$
is constructed from the $k-1$st braid by conjugating by 
$\tau_{2k+1,m+k+1}$,  for $k=1,\dots,m-1$. 

\begin{figure}[htbp]
\begin{center}
\includegraphics[height=1.4in]{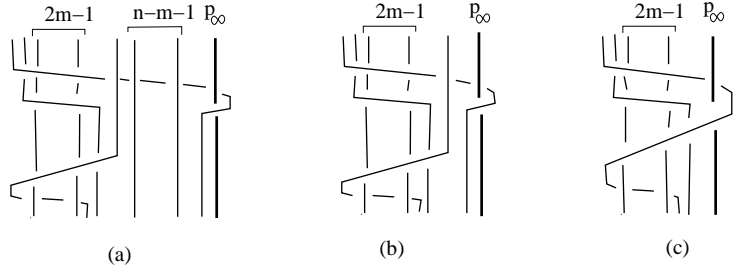}
\caption{After $m-1$ inductive steps: (a) $n \ge m+2$, (b) $n = m+1$, (c) $n=m$ .}
\label{interbraid-fig}
\end{center}
\end{figure}

The resulting braid $\widehat{\sigma_{m,n}}^{(m-1)}$ takes one of three forms: Figure~\ref{interbraid-fig}
(a) shows the general case when $n \ge m+2$, Figure~\ref{interbraid-fig}(b) shows the
case when $n=m$, and Figure~\ref{interbraid-fig}(c) shows the case when $n=m+1$.

\begin{figure}[htbp]
\begin{center}
\includegraphics[height=1.4in]{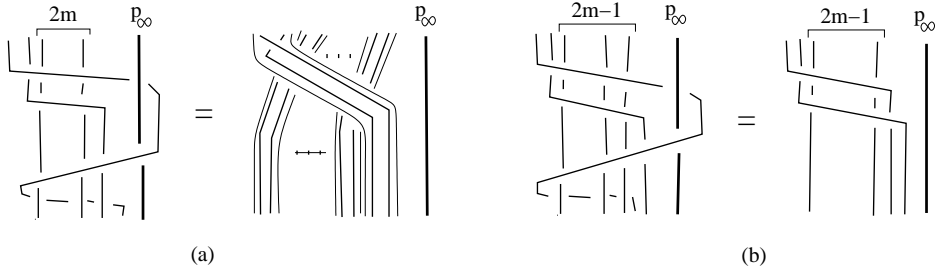}
\caption{Reducible and periodic cases.}
\label{perred-fig}
\end{center}
\end{figure}

\begin{prop}\label{reducible-prop}
When $n=m+1$,
$\sigma_{m,n}$ is a reducible braid.
\end{prop}

\proof  By applying one more conjugation by $\tau_{2m+1,2m+1}$, we obtain the right-hand braid
in  Figure~\ref{perred-fig}(a).  One sees that there is a collection of disjoint disks enclosing pairs of
marked points in $S^2$ whose boundaries are left invariant by $\phi_{\widehat{\sigma_{m,n}}}$.
The claim now follows from Lemma~\ref{boundary-lem}. \qed
\medskip

\begin{prop}\label{periodic-prop}
When $n=m$, $\sigma_{m,n}$ is a periodic braid.
\end{prop}

\proof  Figure~\ref{perred-fig}(b) shows an equivalence of spherical braids.  It is not hard to see
that the right-hand braid is periodic in $\sB(S^2;\sS,\{p_\infty\})$.  The rest follows from
Lemma~\ref{boundary-lem}. \qed
\medskip

\begin{figure}[htbp]
\begin{center}
\includegraphics[width=4.5in]{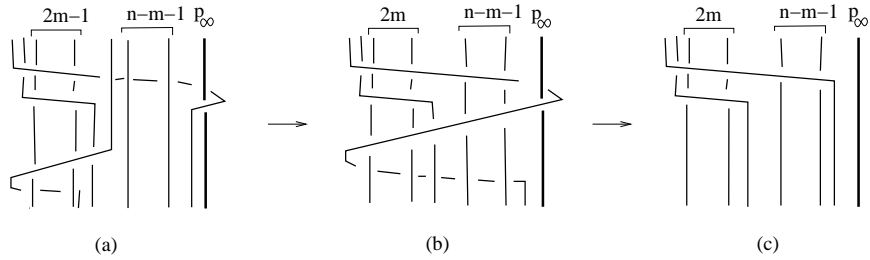}
\caption{General case.}
\label{pA-fig}
\end{center}
\end{figure}

\noindent
The general case, when $n \ge m+2$,  is shown in Figure~\ref{pA-fig}.  The transition from
Figure~\ref{pA-fig}(a) to \ref{pA-fig}(b) is given by doing successive conjugations by
$\tau_{2m+k,2m+k}$ for $k=1,\dots,n-m$.

\begin{figure}[htbp]
\begin{center}
\includegraphics[width=5in]{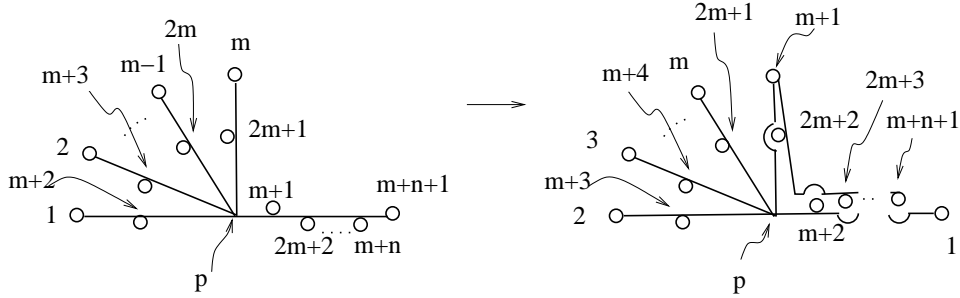}
\caption{Graph map $\h_{m,n}$ for $\phi_{\sigma_{m,n}}$.}
\label{SmnGraph-fig}
\end{center}
\end{figure}

Now consider the graph map $h_{m,n} : H_{m,n} \rightarrow \sF(H_{m,n})$ drawn 
in Figure~\ref{SmnGraph-fig} associated to the preimage of the braid in Figure~\ref{pA-fig}(c)
under the contraction map of Lemma~\ref{boundary-lem}.  
The unusual numbering comes from the left-to-right ordering
of the strands (excluding $p_\infty$) of $\widehat{\sigma_{m,n}}$ 
shown in Figure~\ref{spherical-fig}(b) (cf. end of Section~\ref{TT-section}).  This ordering
proves useful for comparing the transition matrices of $\phi_{\beta_{m,n}}$ and
$\phi_{\sigma_{m,n}}$ (see Section~\ref{Salem-Boyd-section}).

 \begin{prop}\label{SmnGraph-prop} For $n \ge m+2$,
 the graph map $\h_{m,n} : H_{m,n} \rightarrow \sF(H_{m,n})$ is a graph map
 for a conjugate of $\phi_{\sigma_{m,n}}$ satisfying (BH:1) and (BH:2).
 \end{prop}

\proof
One can see that $\h_{m,n}$ is a graph map for $\phi_{\sigma_{m,n}}$ by 
looking at the conjugate form of $\widehat{\sigma_{m,n}}$ given in  Figure~\ref{pA-fig} (c).
The proof that $\h_{m,n}$ satisfies (BH:1) and (BH:2) 
is similar to that of Proposition~\ref{BmnGraph-prop}.\qed
\medskip

\noindent
{\bf Proof of Theorem~\ref{Smn-classification-thm}.}  By Proposition~\ref{Smn-sym-prop}, it
suffices to classify the braids $\sigma_{m,n}$ with $n \ge m \ge 1$.  By Proposition~\ref{reducible-prop},
$\sigma_{m,n}$ is reducible if $n = m+1$, and by Proposition~\ref{periodic-prop}
$\sigma_{m,n}$ is periodic if $n=m$.  In all other cases, Proposition~\ref{SmnGraph-prop} and Proposition~\ref{Smn-sym-prop}  imply
that $\sigma_{m,n}$ is pseudo-Anosov and satisfies
$\lambda(\sigma_{m,n}) = \lambda(\sigma_{n,m})$.\qed
\medskip

\subsection{Train tracks} \label{TT-section}

\begin{figure}[htbp]
\begin{center}
\includegraphics[angle=0,height=1.5in]{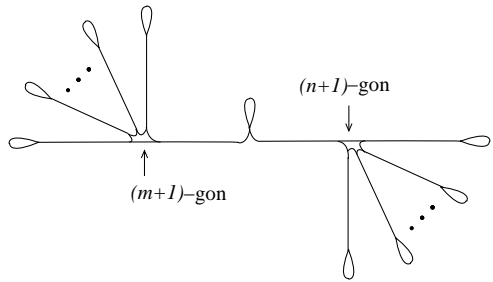}
\caption{Train track for $\phi_{\beta_{m,n}}$.}
\label{Bmn-TT}
\end{center}
\end{figure}

In this section, we find train tracks for $\phi_{\beta_{m,n}}$ and 
$\phi_{\sigma_{m,n}}$, and derive properties of these mapping classes.
The train track for $\phi_{\beta_{m,n}}$ 
associated to $\g_{m,n}$ is given in Figure~\ref{Bmn-TT}.

\begin{lem}\label{Bmn-prongs-lem} The invariant foliations associated to 
the pseudo-Anosov representative of 
 $\phi_{\widehat{\beta_{m,n}}} \in \sM(S^2; \sS,\{p_\infty\})$ have 1-pronged singularities at 
  $p_\infty$ and marked points, and 
 an $(m+1)$-  (respectively, $(n+1)$-) pronged singularity at the fixed point 
$p$ (respectively, fixed point $q$).
\end{lem}

\proof Apply  Lemma~\ref{prongs-lem} to
the train track in Figure~\ref{Bmn-TT}.\qed
\medskip

\begin{figure}[htbp]
\begin{center}
\includegraphics[height=1.5in]{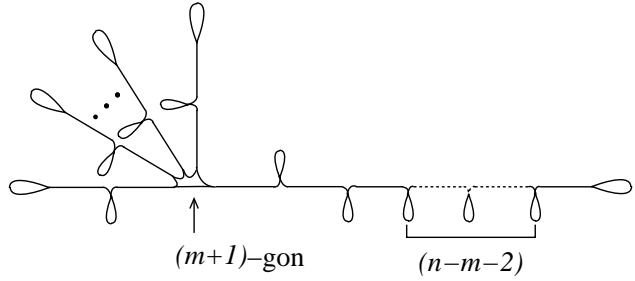}
\caption{Train track for $\phi_{\sigma_{m,n}}$.} 
\label{Smn-TT}
\end{center}
\end{figure}

\noindent
The train track for $\phi_{\sigma_{m,n}}$ 
associated to $\h_{m,n}$ is given in 
Figure~\ref{Smn-TT}.

\begin{lem}\label{Smn-prongs-lem}  For $n \ge m+2$, the invariant 
foliations associated to the pseudo-Anosov representative of
$\phi_{\widehat{\sigma_{m,n}}} \in \sM(S^2; \sS, \{p_{\infty}\})$ have 1-pronged singularities 
at the marked points, an $(m+1)$-pronged
singularity at the fixed point $p$, and an $n$-pronged singularity  at  $p_{\infty}$.
\end{lem}

\proof Apply Lemma~\ref{prongs-lem} to the train track in Figure~\ref{Smn-TT}.\qed
\medskip

By Lemmas \ref{closure2-lem} and \ref{Smn-prongs-lem}, we have the following. 

\begin{cor}\label{Overline-Smn-lem}  For $n \ge m+2$, $\overline{\sigma_{m,n}}$ is
pseudo-Anosov, and
$$
\lambda(\sigma_{m,n}) = \lambda(\overline{\sigma_{m,n}}).
$$
\end{cor}

A pseudo-Anosov map is said to be {\it orientable} if its stable and unstable foliations have only
even order prongs.   

\begin{prop}\label{Bmn-lift-prop} 
Let $m+n=2g$. 
If both $m$ and $n$ are odd, 
there is a pseudo-Anosov element of $\sM_g^0$ whose pseudo-Anosov representative is 
orientable with the same dilatation as $\beta_{m,n}$. 
\end{prop}

\proof  
Let $\Phi_{\widehat{\beta_{m,n}}}$ be the pseudo-Anosov representative of $\phi_{\widehat{\beta_{m,n}}}$, and let
$\Phi'_{\widehat{\beta_{m,n}}}$ be its lift to the double branched covering $F_g$ of $S^2$ branched along
$\sS \cup \{p_\infty\}$.  By Lemma~\ref{Bmn-prongs-lem}, $\Phi_{\widehat{\beta_{m,n}}}'$ 
is a pseudo-Anosov map on $F_g$ with
invariant foliations having two $m+1$- (respectively, $n+1$-) pronged singularities above $p$
(respectively, $q$).\qed

\begin{prop}\label{Smn-lift-prop} 
There is a pseudo-Anosov element of $\sM_g^0$ whose pseudo-Anosov representative is 
orientable with the same dilatation as $\sigma_{m,n}$, 
for each $m,n \ge 1$ with $m+n=2g$ and $|m-n| \ge 2$. 
\end{prop}

\proof  By Proposition~\ref{Smn-sym-prop}, we can assume $n \ge m+2$.
Since $m$ and $n$ have the same parity, Lemma~\ref{Smn-prongs-lem} implies that 
for the invariant foliations associated to the pseudo-Anosov representative $\Phi_{\widehat{\sigma_{m,n}}}$ 
of $\phi_{\widehat{\sigma_{m,n}}}$,  the number of prongs $n_p$ (respectively, $n_\infty$)
at $p$ (respectively $p_\infty$) have opposite parity.  

Let $F_g$ be the branched covering of $S^2$ branched along $\sS$ and either
$p$ if $n_p$ is odd, or $p_\infty$ if $n_\infty$ is odd.  Let $\Phi_{\widehat{\sigma_{m,n}}}'$ be the
lift of $\Phi_{\widehat{\sigma_{m,n}}}$  to $F_g$.    Then $\Phi_{\widehat{\sigma_{m,n}}}'$ is pseudo-Anosov
with dilatation equal that of $\Phi_{\widehat{\sigma_{m,n}}}$.   Furthermore, by our choice of
branch points, the invariant foliations have only even order prongs, and hence are
orientable.\qed
\medskip

We conclude this section by relating the graph maps of $\phi_{\beta_{m,n}}$ and
$\phi_{\sigma_{m,n}}$ in a way that is compatible with the conjugations used in
Section~\ref{symmetry-section}.  
\begin{figure}[htbp]
\begin{center}
\includegraphics[width=5in]{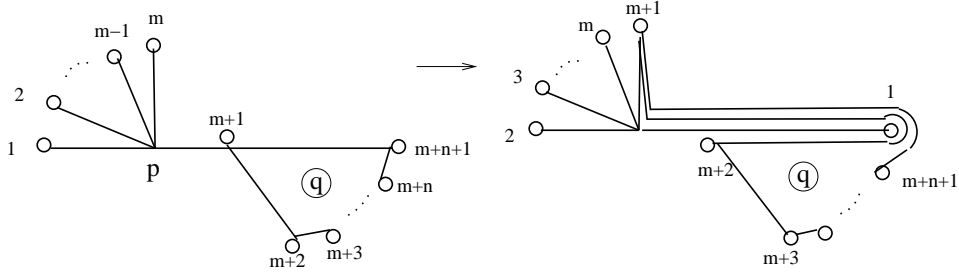}
\caption{Graph map $g_{m,n}'$ for $\phi_{m,n}$.}
\label{BmnGraphNew-fig}
\end{center}
\end{figure}

Let $\Phi_{\beta_{m,n}}$ be the pseudo-Anosov representative
 of $\phi_{\beta_{m,n}}$.  Since
$q$ is a fixed point for $\Phi_{\beta_{m,n}}$,
$\Phi_{\beta_{m,n}}$ defines a mapping class
 $\phi_{m,n} = [\Phi_{\beta_{m,n}}]$ in 
$\sM(S^2;\sS,\{q\}, \{p_\infty\})$.
Identify $g_{m,n}$ with the graph map on $S^2$
for $\phi_{m,n}$ obtained by pushed forward
by the contraction map in Lemma~\ref{boundary-lem}.
Let 
$$\g_{m,n}' : G_{m,n}' \rightarrow G_{m,n}'$$ be the graph map 
obtained from $g_{m,n}$ after 
 puncturing the sphere at the fixed point $q$ for $\g_{m,n}$ 
as in  Figure~\ref{BmnGraphNew-fig}.     

\begin{figure}[htbp]
\begin{center}
\includegraphics[width=5in]{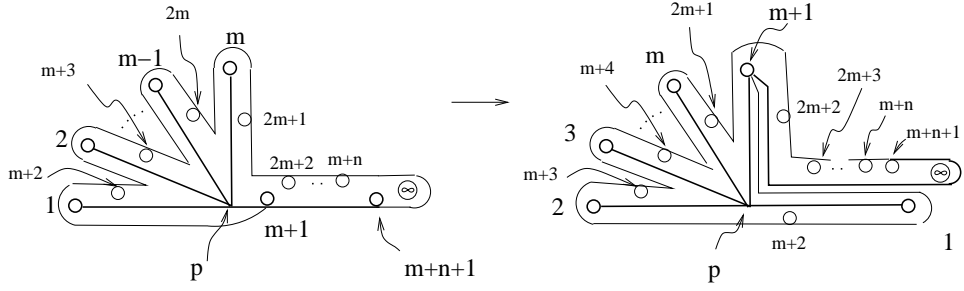}
\caption{Exchanging the roles of $q$ and $p_\infty$ for $g_{m,n}'$.}
\label{wrap-fig}
\end{center}
\end{figure}
\begin{figure}[htbp]
\begin{center}
\includegraphics[width=5in]{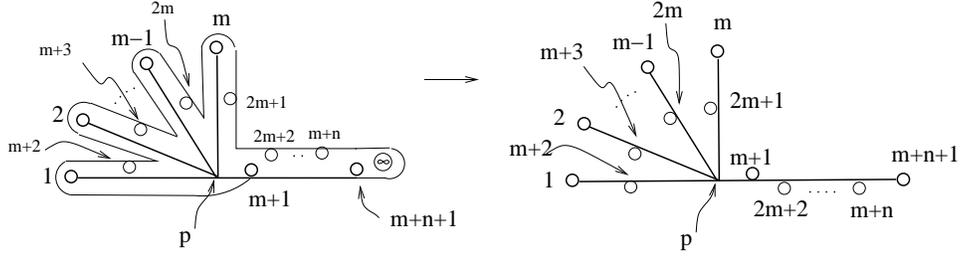}
\caption{Identifying edges on $G_{m,n}'$.}
\label{identify-fig}
\end{center}
\end{figure}

 Exchanging the roles of $q$ and $p_\infty$ (i.e., bringing
$p_\infty$ into the visual plane) yields the
 graph map shown in  Figure~\ref{wrap-fig}, which is equivalent
to $g_{m,n}'$.
Now remove the marking at $p_\infty$ and consider the graph map
\begin{eqnarray}\label{graphmap-map}
\f_{m,n} : G_{m,n}' \rightarrow \mathcal{F}(H_{m,n})
\end{eqnarray}
obtained by a natural identification of edges shown in  Figure~\ref{wrap-fig}.
Figure~\ref{identify-fig} shows the projection map $\pi$ applied to the image of $\f_{m,n}$.
 The graph map $\h_{m,n} : H_{m,n} \rightarrow \sF(H_{m,n})$ in Figure~\ref{SmnGraph-fig}
is the one induced by pushing forward
$\g_{m,n}'$  by the map $\f_{m,n}$.

\subsection{Characteristic equations for dilatation}\label{dilatations-section}

\begin{figure}[htbp]
\begin{center}
\includegraphics[height=0.9in]{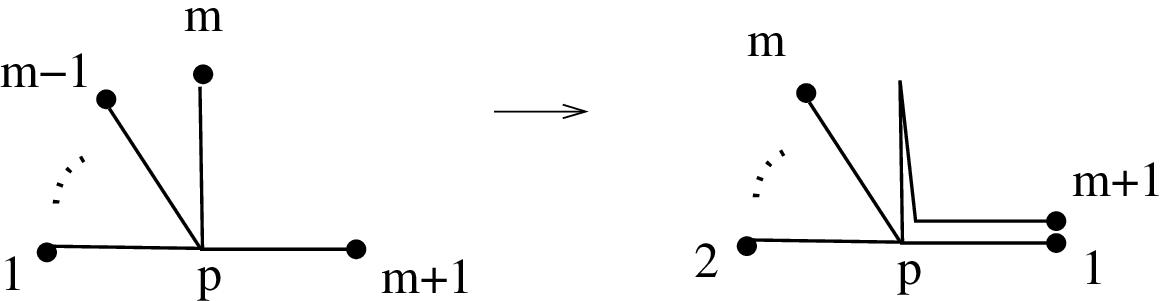}
\caption{The graph map $\r_m : L_m \rightarrow \sF(L_m)$.}
\label{RmnGraph-fig}
\end{center}
\end{figure}

Consider the graph map $\r_m : \Gamma_m \rightarrow \sF(\Gamma_m)$, 
shown in Figure~\ref{RmnGraph-fig}.
As seen in Figures~\ref{BmnGraph-fig} and \ref{SmnGraph-fig}, the graph maps for 
 $\phi_{\beta_{m,n}}$ and $\phi_{\sigma_{m,n}}$ ``contain" $\r_m$ as the 
action on a subgraph.

The transition matrix for $\r_m$ has the following form with respect to the basis of
edges $e(p,1), \dots, e(p,m+1)$: 
$$
\sR_m = \left[
\begin{array}{cccccc}
0 & 1  &  0 & \dots&0 & 0\\
0 & 0 & 1 & \dots &0& 0\\
\dots\\
0 &0 &0& \dots & 1 & 0\\
0 &0 & 0 & \dots & 0& 2\\
1 &0 &0& \dots & 0 & 1\\
\end{array}
\right ].
$$
\medskip

\noindent
The characteristic polynomial for $\sR_m$ is $R_m(t) = t^m(t-1) - 2$.
As we will see in the proof of Theorem~\ref{CharEq-thm} given below,
the appearance of $\sR_m$ 
 within the transition matrices of $\phi_{\beta_{m,n}}$ and
 $\phi_{\sigma_{m,n}}$ gives rise to a similar form for their characteristic equations.
 Given a polynomial $f(x)$ of degree $d$, let $f_*(t) = t^df(1/t)$.
 
\begin{thm}\label{CharEq-thm}  The dilatations of $\beta_{m,n}$ and $\sigma_{m,n}$
are as follows:
\begin{description}
\item{(1)} For $m,n \ge 1$, the dilatation $\lambda(\beta_{m,n})$ is
the largest root of the polynomial
$$
T_{m,n}(t) = t^{n+1} R_m(t) + (R_m)_*(t).
$$
\item{(2)} For $m,n\ge 1$ and $|m-n| \ge 2$, the dilatation $\lambda(\sigma_{m,n})$ is
the largest root of the polynomial
$$
S_{m,n}(t) = t^{n+1} R_m(t) - (R_m)_*(t).
$$
\end{description}
\end{thm}

\proof
In the discussion at the end of Section~\ref{TT-section},
we saw that, up to isotopy, 
$\g_{m,n}'$ is a graph map  for a pseudo-Anosov representative of
$\phi_{\beta_{m,n}}$, and by Proposition~\ref{SmnGraph-prop}
$\h_{m,n}$ is a graph map for a pseudo-Anosov representative of
$\phi_{\sigma_{m,n}}$.   
Thus, to find characteristic polynomials for $\lambda(\beta_{m,n})$ and $\lambda(\sigma_{m,n})$, 
it is enough to compute the transition matrices for $\g_{m,n}'$ and $\h_{m,n}$.

Consider the  basis for $V(G_{m,n}')$:
 \begin{eqnarray*}
v_k &=& e(p,k),\quad k=1,\dots,m, \\
v_{m+1} &=& e(p,m+n+1), \\
v_{m+1+k} &=& e(m+k,m+k+1), \quad k=1,\dots,n, \\
v_{m+n+2} &=& e(p,m+1).
\end{eqnarray*}
The corresponding transition matrix for $g_{m,n}'$ is given by 
$$
\sT'_{m,n} = 
\left [
\begin{array}{cccccc|cccccc}
0 & 1 & 0 & \dots & 0 & 0 &0&0&\dots&&0\\
0 & 0 & 1 & \dots & 0 & 0&0&0& \dots &&0\\
&&&\dots&&&&&\dots&&\\
0 &0 &0&\dots&1&0&0&0&\dots&&0\\
0 & 0 & 0 & \dots & 0 & 2 &1&0& \dots&&(1)_b\\
1 & 0 & 0 & \dots & 0 & 1 &2 & 0 &\dots &&0\\
\hline
0 &0& \dots &&0&0 &0&1&0&\dots &0\\
0 &0& \dots &&0&0&0&0&1&\dots &0\\
&&\dots&&&&&&\dots&\\
0 &0& \dots &&0&0&0&0&\dots &1&0\\
0 &0& \dots &&0&(1)_{ab}&0&0&\dots &0&0\\
0&0&\dots &&0&0&(-1)_b&0&\dots &0&(0)_a\\
\end{array}
\right ]
$$

We will show that the characteristic polynomial for  $\sT'_{m,n}$ is given by
$$
T_{m,n}(t) = t^n R_m(t) + (R_m)_*(t).
$$

The upper left  block matrix of $\sT'_{m,n}$ corresponding to the vectors
$v_1,\dots,v_{m+1}$
is identical to $\mathcal{R}_m$.  Multiplying the characteristic polynomials of the upper
left and lower right
diagonal blocks
gives $t^{n+1}R_m$.  The rest of the characteristic polynomial has two nonzero
summands.  One corresponds to the matrix entries marked $a$, and is given by
$$
t(-1)^{n+1}
\left |
\begin{array} {ccccc}
-1&0&\cdots&0&0\\
t&-1&\cdots&0&0\\
&&\cdots&&\\
0&0&\cdots&t&-1\\
\end{array}
\right |_{(n-1)\times (n-1)}
\left |
\begin{array}{cccccc}
t&-1&\cdots&0&0\\
0&t&\cdots&0&0\\
&&\cdots&&\\
0&0&\cdots&t&-1\\
-1&0&\cdots&0&-2
\end{array}
\right |_{(m+1)\times (m+1)},
$$
which yields $-t(2t^m +1)$.  The other summand corresponds to the matrix entries marked $b$ 
and is given by 
$$
(-1)^{n+1}
\left |
\begin{array} {ccccc}
-1&0&\cdots&0&0\\
t&-1&\cdots&0&0\\
&&\cdots&&\\
0&0&\cdots&t&-1\\
\end{array}
\right |_{(n-2)\times (n-2)}
\left |
\begin{array}{cccccc}
t&-1&0&\cdots&0&0\\
0&t&-1&\cdots&0&0\\
&&&\cdots&&\\
0&0&0&\cdots&t&-1\\
-1&0&0&\cdots&0&0\\
\end{array}
\right |_{m \times m},
$$
which yields $1$.
This completes the proof of (1). 

Let $\sS_{m,n}$ be the transition matrix for $\h_{m,n}: H_{m,n} \rightarrow \sF(H_{m,n}) $.  
We will pull back $\sS_{m,n}$ to an invertible linear transformation on $V(G_{m,n}')$ using the map
$\f_{m,n}$ given in (\ref{graphmap-map}).   
Let $\h_{m,n}' (v_i)$ be the image of $\g_{m,n}'(v_i)$ under the identification defined
by $\f_{m,n}$.  
Then the transition matrix $\sS_{m,n}'$ for  $\h_{m,n}': G_{m,n}' \rightarrow F(H_{m,n})$  
is the same as $\sT_{m,n}'$ except at the vector $v_{m+n+1}$.  As can be seen
in Figure~\ref{wrap-fig}, we have
$$
\sS_{m,n}'(v_{m+n+1}) = \sT_{m,n}'(v_{m+1}) - 2v_{m+n+1}
$$ 
Thus, $\sS_{m,n}'$ differs from $\sT_{m,n}'$ only by changing the entry labeled by both
$a$ and $b$ from $1$ to $-1$.

Recall that  the sign of the
entry marked both $a$ and $b$ in $\sT_{m,n}$ determines the sign of
in front of $(R_{m})_*$.  Since this sign is the only difference between $\sS_{m,n}'$
and $\sT_{m,n}$, the characteristic polynomial for $\sS'_{m,n}$ is given
by
$$
S_{m,n}(t) = t^nR_m(t) - (R_m)_*(t).
$$

To finish the proof of (2) we have left to check that $\lambda(\sigma_{m,n})$ is 
the largest root of $S_{m,n}$.
Thus (2) follows if we can show that the extra eigenvalue of $\sS_{m,n}'$ has 
absolute value $1$.  We will show below that the
extra eigenvalue equals $1$.

From Figure~\ref{identify-fig}, we see that the kernel of the linear map induced by
$\f_{m,n}$ is spanned by
$$
w = 2(v_1 + \cdots + v_m) + v_{m+1} - (v_{m+2} + \cdots + v_{m+n+1})  + v_{m+n+2}.
$$
Under the map $\h_{m,n}'$, we have
\begin{eqnarray*}
2(v_1+ \cdots v_m) &\mapsto& 2(v_{m+1} + v_1 + \cdots + v_{m-1}),\\
v_{m+1} &\mapsto& 2v_m + v_{m+1} - v_{m+n+1},\\
v_{m+2} + \cdots +v_{m+n+1} &\mapsto& v_{m} + 2 v_{m+1}  + 
v_{m+2} + \cdots + v_{m+n} - v_{m+n+2},\\
v_{m+n+2} &\mapsto& v_{m},
\end{eqnarray*}
and hence, $\h_{m,n}'(w) =  w$.   Thus, the characteristic polynomial for $\sS_{m,n}'$
differs from that for $\sS_{m,n}$ by a factor of $(t-1)$.
\qed

\begin{rem} Hiroyuki Minakawa independently discovered the
pseudo-Anosov maps on $F_g$ constructed in the proof of Proposition~\ref{Smn-lift-prop}, for the
case when  $(m,n) = (g-1,g+1)$, using a beautiful new method for constructing
orientable pseudo-Anosov maps on genus $g$ surfaces.  
He also directly computes their dilatation using different techniques from ours.
\end{rem}

\subsection{Dilatations and Salem-Boyd sequences}\label{Salem-Boyd-section}
In this section, we apply properties of Salem-Boyd sequences to find the least
dilatations among $ \lambda(\sigma_{m,n})$ and $\lambda(\beta_{m,n})$ 
for $m+n = 2g$ fixing $g$.
We also give bounds on these dilatations. 

Given a polynomial $f(t)$ of degree $d$, the {\it reciprocal} of $f(t)$ is $f_*(t) = t^df(1/t)$.   
The polynomial $f$ satisfying $f=f_*$ (respectively,  $f=-f_*$) 
 is a {\it reciprocal polynomial} (respectively, {\it anti-reciprocal polynomial}).
For a monic integer polynomial $P(t)$ of degree $d$, 
the sequence 
$$
Q^\pm_n(t) = t^n P(t) \pm P_*(t)
$$
is called the  {\it Salem-Boyd sequence} associated to $P$.   

\begin{thm}\label{Salem-Boyd-theorem}
Let $Q_n$ be a Salem-Boyd sequence associated to $P$.
Then $Q_n$ is a reciprocal or an anti-reciprocal polynomial, and the set of roots of
$Q_n$ outside the unit circle converge to those of $P$ as $n$ goes
to infinity.  
\end{thm}
Theorem~\ref{Salem-Boyd-theorem} is a consequence of Rouch\'e's Theorem
applied to the sum
$$
\frac{P(t)}{t^d} \pm \frac{P_*(t)}{t^{n+d}}
$$
considered as a holomorphic function on the Riemann sphere minus the unit disk.

For a  monic integer polynomial $f(t)$, let $N(f)$ be the number of roots of $f$ outside
the unit circle, $\lambda(f)$ the maximum norm of roots of $f$, and $M(f)$ the
product of the norms of roots outside the unit circle, 
which is called the {\it  Mahler measure} of $f$. 

\begin{cor}\label{Salem-Boyd-cor}  Let $Q_n$ be a Salem-Boyd sequence
associated to $P$.  Then 
\begin{description}
\item{(1)} $\lim_{n\rightarrow \infty} M(Q_n) = M(P)$, and
\item{(2)} $\lim_{n\rightarrow \infty} \lambda(Q_n) = \lambda(P)$.
\end{description}
\end{cor}

Any algebraic integer on the unit circle has a (anti-)reciprocal minimal polynomial.  Suppose
$$
P(t) = P_0(t) R(t),
$$
where $R$ is (anti-)reciprocal and $P_0$ has no roots on the unit circle.
Then 
$$
Q_n(t) = R(t) (t^nP_0(t) \pm (P_0)_*(t)).
$$
We have thus shown the following.

\begin{lem}
All roots of $P$ on the unit circle are also roots of $Q_n$ for all $n$.
\end{lem}

The following theorem can be proved by first restricting to the case when $P$
has no roots on the unit circle, and then by defining a natural deformation of
the roots of $P(t)$ to those of $Q_n(t)$, which don't cross the unit
circle (see \cite{Boyd77}).

\begin{thm}\label{number-of-roots-theorem}
Let $Q_n$ be a Salem-Boyd sequence associated to $P$.  Then 
$$
N(Q_n) \leq N(P)
$$
for all $n$.
\end{thm}

We now apply the above results to the Salem-Boyd sequences $S_{m,n}$ and
$T_{m,n}$ of Theorem~\ref{CharEq-thm}.  To do this we first study the polynomials
$R_m$.

\begin{lem}\label{MRm-lem}  For all $m$, $M(R_m) = 2$.
\end{lem}

\proof  For $|t| < 1$ we have
$$
|t^m(t-1)| < 2,
$$
and hence $R_m$ has no roots strictly within the unit circle.
Therefore, the Mahler measure of $R_m$ must equal the absolute value of the constant
coefficient, namely 2. \qed

\begin{cor}\label{asymp-cor} Fixing $m$ and letting $n$ increase,
the Mahler measures of $T_{m,n}$ and $S_{m,n}$ converge to 2.
\end{cor}

\proof Apply Corollary~\ref{Salem-Boyd-cor}.\qed
\medskip

\begin{lem}\label{Rm-lem}  The polynomial $R_m$ has one real root
outside the unit circle.  This root  is simple and greater than 1.
\end{lem}

\proof Taking the derivative
$$
R_m'(t) = (m+1)t^m - mt^{m-1} 
$$ 
we see that $R_m$ is increasing for $t > \frac{m}{m+1}$, and hence also for $t \ge 1$.
Since $R_m(1) = -2 < 0$ and
$R_m(2) > 0$, it follows that $R_m$ has a simple root $\mu_m$ with
$1 < \mu_m < 2$.  Similarly we can show that, for $t < 0$, $R_m$ has no roots for
$m$ even, and one root if $m$ is odd.  In the odd case, $R_m(-1) = 0$,
so $R_m$ has no real roots strictly less than $-1$. \qed

\begin{lem}\label{Rm2-lem}  The sequence
$\lambda(R_m)$ converges monotonically to 1 from above.
\end{lem}

\proof
 Since $M(R_m) = 2$, we know that $\mu_m = \lambda(R_m) > 1$.  
Take any $\epsilon > 0$. 
Let $D_\epsilon$ be the disk of radius $1+\epsilon$ around the origin in the complex
plane. Let $g(t) = \frac{t-1}{t}$ and $h_m(t) = \frac{-2}{t^{m+1}}$.  Then for large enough $m$, we 
have
$$
|g(t)| = \left | \frac{t-1}{t} \right | > \left | \frac{2}{t^{m+1}} \right | = |h_m(t)|
$$
for all $t$ on the boundary of $D_\epsilon$,
and $g(t)$ and $h_m(t)$ are holomorphic on the complement of $D_\epsilon$ in 
the Riemann sphere.  By Rouch\'e's theorem, 
$g(t)$, $g(t) + h_m(t)$, and hence $R_m(t)$
have the same number of roots outside $D_\epsilon$,
which is zero.

To show monotonicity consider $R_m(\mu_{m+1})$.  
Since $\mu_{m+1}$ satisfies
$$
(\mu_{m+1})^{m+1}(t-1) -2 = 0,
$$
we have
\begin{eqnarray*}
R_m(\mu_{m+1}) &=&( \mu_{m+1})^m (t-1) -2\\
&=& ((\mu_{m+1})^m - (\mu_{m+1})^{m+1}) (t-1)\\
&<& 0.
\end{eqnarray*}
Since $R_m(t)$ is an increasing function for $t> 1$, we
conclude that $\mu_{m+1}< \mu_m$.
\qed

\begin{cor}
For fixed $m$, the sequences $\lambda(\beta_{m,n})$ and $\lambda(\sigma_{m,n})$
converge to $\lambda(R_m)$ as sequences in $n$.  Furthermore, we can make
$\lambda(\beta_{m,n})$ and $\lambda(\sigma_{m,n})$ arbitrarily small by taking
$m$ and $n$ large enough.  \end{cor}

We now determine the monotonicity of $\lambda(\beta_{m,n})$ and
$\lambda(\sigma_{m,n})$ for fixed $m$.

\begin{prop}\label{Perron-prop}
For fixed $m$,  the dilatations $\lambda(\beta_{m,n})$ are strictly monotone
decreasing, and for $n \ge m+2$, the dilatations $\lambda(\sigma_{m,n})$ are 
strictly monotone increasing.
\end{prop}

\proof  
Consider 
$$
f(t) = (R_m)_*(t) = -2t^{m+1} - t + 1.
$$
Then, for $t > 0$,
$$
f'(t) = -2(m+1)t^m - 1 < 0.
$$
Also $f(1) = -2 < 0$.
Since
$b_{m,n} = \lambda(\beta_{m,n}) > 1$, and for $n \ge m+2$, $s_{m,n} = \lambda(\sigma_{m,n}) > 1$,
it follows that 
$(R_m)_*(b_{m,n})$ and $(R_m)_*(s_{m,n})$ are both negative.
 
We have
$$
0 = T_{m,n}(b_{m,n}) = (b_{m,n})^{n+1} R_m(b_{m,n}) + (R_m)_*(b_{m,n}),
$$
and
$$
0 = S_{m,n}(s_{m,n}) = (s_{m,n})^{n+1} R_m(s_{m,n}) - (R_m)_*(s_{m,n}),
$$
which imply that
$$
R_m(b_{m,n})  >  0\quad \mbox{and} \quad R_m(s_{m,n}) < 0.
$$
Since $R_m$ is increasing for $t > 1$, we have 
$b_{m,n} > \mu_m$ and $s_{m,n} < \mu_m$.

Plug $b_{m,n}$ into $T_{m,n-1}$, and subtract $T_{m,n}(b_{m,n}) = 0$:
\begin{eqnarray*}
T_{m,n-1}(b_{m,n}) &=&(b_{m,n})^{n-1}R_m(b_{m,n}) + (R_m)_*(b_{m,n})\\
&=&((b_{m,n})^{n-1} - (b_{m,n})^n)R_m(b_{m,n})\\
&<& 0.
\end{eqnarray*}
Since $b_{m,n-1}$ is the largest real root of $T_{m,n+1}$, we
have $b_{m,n} < b_{m,n-1}$.
We can show that $s_{m,n} < s_{m,n+1}$ for $n \ge m+2$ in a similar way, 
by adding the formula for $S_{m,n}(s_{m,n})$ to $S_{m,n+1}(s_{m,n})$.
\qed

\begin{cor} 
\label{inequality-cor}
We have the inequality 
$$
\lambda(\beta_{m,n})> \lambda(\sigma_{m,n})
$$
for all $m,n$ with $|m-n| \ge 2$. 
\end{cor}

We now fix $2g=m+n$, and show that among the braids $\beta_{m,n}$ and
$\sigma_{m,n}$ for $g \ge 2$, $\sigma_{g-1,g+1}$ has least dilatation.

\begin{prop}\label{min-prop}  The braids $\beta_{m,n}$ and
$\sigma_{m,n}$ satisfy
\begin{eqnarray*}
\lambda(\beta_{m,m}) &<& \lambda(\beta_{m-k,m+k}), \\
\lambda(\beta_{m,m+1}) &<& \lambda(\beta_{m-k,m+k+1}) 
\end{eqnarray*}
for $k=1,\dots,m-1$, and 
\begin{eqnarray*}
\lambda(\sigma_{m-1,m+1}) &<& \lambda(\sigma_{m-k,m+k}), \\
\lambda(\sigma_{m-1,m+2}) &<& \lambda(\sigma_{m-k,m+k+1})
\end{eqnarray*}
for $k=2,\dots,m-1$.
\end{prop}

\proof  Let $\lambda= \lambda(\beta_{m,m})$.  Then plugging $\lambda$ into
$T_{m-k,m+k}$ gives
\begin{eqnarray*}
T_{m-k,m+k}(\lambda) &=& \lambda^{m+k+1}(\lambda^{m-k} (\lambda - 1) - 2) - 2 \lambda^{m-k+1} - \lambda + 1\\
&=& \lambda^{2m+2} - \lambda^{2m+1} - 2\lambda^{m+k+1} - 2\lambda^{m-k+1} - \lambda +1. 
\end{eqnarray*}
Subtracting 
$$
0 = T_{m,m}(\lambda) = \lambda^{2m+2} - \lambda^{2m+1} - 4\lambda^{m+1} - \lambda + 1, 
$$
we obtain
\begin{eqnarray*}
T_{m-k,m+k}(\lambda) &=& 4 \lambda^{m+1} - 2\lambda^{m+k+1} - 2\lambda^{m-k+1}\\
&=& -2\lambda^{m-k+1}(\lambda^k-1)^2\\
&<& 0.
\end{eqnarray*}
Since $\lambda(\beta_{m-k,m+k})$ is the largest real root of $T_{m-k,m+k}$, we have
$\lambda(\beta_{m,m})<\lambda(\beta_{m-k,m+k})$.
The other inequalities are proved similarly.
\qed

\begin{prop}\label{min-Smn-Bmm-prop}
For $m \ge 2$, the dilatations of $\beta_{m,n}$ and $\sigma_{m,n}$ satisfy the inequalities
$$
\lambda(\beta_{m,m}) > \lambda(\sigma_{m-1,m+1})
$$
and
$$
\lambda(\beta_{m,m+1}) \ge \lambda(\sigma_{m-1,m+2})
$$
with equality if and only if $m=2$.
\end{prop}

\proof Let $\lambda=\lambda(\sigma_{m-1,m+1})$.  Then
$$
T_{m,m}(\lambda) = \lambda^{2m+2} - \lambda^{2m+1} - 4\lambda^{m+1} - \lambda + 1.
$$
Plugging in the identity
$$
0 = S_{m-1,m+1}(\lambda) = \lambda^{2m+2}-\lambda^{2m+1} - 2\lambda^{m+2} + 2\lambda^m + \lambda - 1,
$$
and subtracting this from $T_{m,m}(\lambda)$ we have
\begin{eqnarray*}
T_{m,m}(\lambda) &=& 2\lambda^{m+2} - 4\lambda^{m+1} - 2\lambda^m - 2\lambda + 2\\
&=& 2\lambda^m(\lambda^2 - 2\lambda + 1) + 2(1-\lambda).
\end{eqnarray*}
The roots of $t^2 - 2t +1$ are $1 \pm \sqrt{2}$ by the quadratic formula.  Since
$$
1 - \sqrt{2} < 1 < \lambda < 2 < 1 + \sqrt{2},
$$
$\lambda^2-2\lambda+1$ and $1-\lambda$ are both negative, hence $T_{m,m}(\lambda) < 0$.
Since $\lambda(\beta_{m,m})$ is the largest real root of $T_{m,m}(t)$, it follows that
$\lambda(\sigma_{m-1,m+1}) = \lambda < \lambda(\beta_{m,m})$.

For the second inequality we plug in $\lambda = \lambda(\sigma_{m-1,m+2})$ 
into $T_{m,m+1}$.  This gives
$$
T_{m,m+1}(\lambda)=2\lambda^m(\lambda^3 - \lambda^2 - \lambda - 1) - \lambda  - 1.
$$
Thus,  $\lambda^3 - \lambda^2 - \lambda - 1 < 0$ would imply $T_{m,m+1} (\lambda) < 0$.  
The polynomial $g(x) = t^3 - t^2 - t -1$
has one real root ($\approx 1.83929$) and is increasing for $t >1$.
Since $\lambda(R_m)$ is decreasing with $m$
and $\lambda < \lambda(R_2) \approx 1.69562 < 1.8$, we 
see that $T_{m,m+1}(\lambda) < 0$ for $m \ge 3$.   For the remaining case, we check that
$T_{2,3}  = S_{1,4}$. 
\qed

\begin{cor}\label{min-braid-cor}  The least dilatation among $\sigma_{m,n}$ and $\beta_{m,n}$ 
for $m+n=2g$ is given by $\lambda(\sigma_{g-1,g+1})$.
\end{cor}

By Corollary~\ref{Salem-Boyd-cor}, 
Lemma~\ref{Rm2-lem} and Proposition~\ref{Perron-prop}, 
 the dilatations $\lambda(\sigma_{m,n})$ for $m+2 \leq n$
converge to 1 as $m,n$ approach infinity.  We prove the following stronger statement,
which implies Theorem~\ref{asymp-thm}.

\begin{prop}\label{min-Smn-prop}
The dilatation $\lambda_g = \lambda(\sigma_{g-1,g+1})$ satisfies
$$
\frac{\log(2+\sqrt{3})}{g+1} < \log (\lambda_g) < \frac{ \log (2 +\sqrt{3})}{g}.
$$
\end{prop}

\proof Using Theorem~\ref{CharEq-thm}, we see that $\lambda = \lambda_g$
satisfies
\begin{eqnarray}\label{dil-eqn}
0&=&\lambda^{2g+1} - 2\lambda^{g+1} - 2\lambda^g + 1\\
&=& \lambda(\lambda^g)^2 - 2(\lambda+1)\lambda^g + 1. \nonumber
\end{eqnarray}
Since $\lambda$ is the largest real solution, the quadratic formula gives
\begin{eqnarray*}
\lambda^g &=& \frac{2(\lambda+1) + \sqrt{4(\lambda+1)^2 - 4\lambda}}{2\lambda}\\
&=& \frac{\lambda+1 + \sqrt{\lambda^2 +  \lambda + 1}}{\lambda}. 
\end{eqnarray*}
It follows that
\begin{eqnarray}\label{dilatation-eqn}
\lambda^{g+1} &=& \lambda+1 + \sqrt{\lambda^2 + \lambda+1}.
\end{eqnarray}
Since $2 > \lambda > 1$ for all $g$, (\ref{dilatation-eqn}) implies
$$
2 + \sqrt{3} < \lambda^{g+1} < 3+\sqrt{7}.
$$

We improve the upper bound using an argument conveyed to us by Hiroyuki Minakawa.
Rewrite the equation (\ref{dil-eqn}) as follows
\begin{eqnarray*}
0 &=&\lambda^{2g+1} + \lambda^{2g} - \lambda^{2g} - 2(\lambda+1)\lambda^g + 1\\
&=&\lambda^{2g}(\lambda+1) - (\lambda^{2g} - 1) - 2(\lambda+1)\lambda^g. 
\end{eqnarray*}
Factoring out  $(\lambda+1)$ gives
$$
0 = \lambda^{2g} - \frac{\lambda^{2g}-1}{\lambda+1} - 2 \lambda^g.
$$
Since $\lambda > 1$, we have
$$
\frac{\lambda^{2g} -1}{\lambda+1} < \frac{1}{2}(\lambda^{2g}-1).
$$
This implies the inequality
$$
x^{2g} - \frac{x^{2g}-1}{x+1} - 2x^g > \frac{1}{2}(x^{2g} - 4x^g + 1) =: p(x)
$$
for $x$ near $\lambda$.  Thus, $p(x)$ has a real root $\mu$ larger than $\lambda$.
Using the quadratic formula again, we see that 
$$
\mu^g = 2 + \sqrt{3},
$$
and hence
$$
\lambda^g < \mu^g = 2+\sqrt{3}.
$$
\qed

\section{Further discussion and questions}\label{discussion-section}

As stated in the introduction,
the braids with smallest dilatations for 3 and 4 strands are $\beta_{1,1}$ and $\beta_{1,2}$,
respectively.  By Proposition~\ref{min-prop} and Proposition~\ref{min-Smn-Bmm-prop},
for $s \ge 5$, the minimal dilatations by our construction
come from $\sigma_{m-1,m+1}$, when $s = 2m+1$; and $\sigma_{m-1,m+2}$, when
$s=2m+2$.   For $s$ even, there are examples of braids with smaller dilatation
than that of $\sigma_{m-1,m+2}$, but for $s$ odd, we know of no such examples.

For $s$ odd, we have $s = 2g+1$, and 
$$
\Sigma(\sB(D,2g+1)) \subset \Sigma(\sM_g^0).
$$
Thus, Penner's lower bound for elements of $\Sigma(\sM_g^0)$ extend to 
$\Sigma(\sB(D,2g+1))$ and we have
$$
\delta(\sB(D,2g+1)) \ge \delta(\sM_g^0) \ge  \frac{\log(2)}{12 g - 12}.
$$
For $g =2$, Zhirov \cite{Zhirov95} shows that if $\phi \in \sM_2^0$ is pseudo-Anosov
with orientable
invariant foliations, then $\lambda(\phi)$ is bounded below by the largest root $\lambda_0$ of
$$
x^4 -x^3 - x^2 -x+1.
$$
For $s=5$, the braid $\sigma_{1,3}$ is pseudo-Anosov, and its lift to a genus
$g=2$ surface is orientable.  Our formula shows that the dilatation of $\sigma_{1,3}$
is the largest root of Zhirov's equation, and hence $\sigma_{1,3}$ achieves the least dilatation
of a genus $2$ orientable pseudo-Anosov map.   This yields the following weaker version
of Ham and Song's result \cite{HamSong05}, which doesn't assume any conditions on the
combinatorics of the invariant foliations.

\begin{cor} The braid $\sigma_{1,3}$ is pseudo-Anosov with smallest dilatation among
all braids on $5$ strands whose invariant folliations are even-pronged at all interior singularities.
\end{cor}

We discuss the following general
 question and related work on the forcing relation in Section~\ref{forcing-section}.
\begin{ques} 
Is there a pseudo-Anosov braid $\beta \in \sB(D, 2g+1)$ such that 
$$
\lambda(\beta)< \lambda(\sigma_{g-1,g+1})\quad ?
$$
\end{ques}

Let $\sK_g^s \subset \sM_g^s$ be the subset of mapping classes that arise as the monodromy of
a fibered link  $(K,F)$ in $S^3$, where the fiber $F$ has genus $g$ and the link $K$ has
$s$ components.
\begin{ques}\label{fiberedlink-ques} Is there a strict inequality
$$
\delta(\sM_g^s) < \delta(\sK_g^s) \quad ?
$$
\end{ques}
In Section~\ref{fiberedQ-section}, we briefly discuss what is known about bounds on dilatations of
pseudo-Anosov monodromies of fibered links, and show how the braids $\beta_{m,n}$ arise in this
class.

\subsection{The forcing relation on the braid types}\label{forcing-section}

The existence of periodic orbits of dynamical systems can imply the existence
of other periodic orbits.
Continuous maps of the interval give typical examples for such phenomena \cite{Sha}.
In \cite{Boy92}, Boyland introduces the notion of {\it braid types}, and defines
a relation on the set of braid types to
study an analogous phenomena  in the $2$-dimensional case.
Let $\sS \subset \mbox{int}(D)$ be a set of $s$ marked points.
There is an isomorphism
$$
\mathcal{B}(D;\sS)/Z(\mathcal{B}(D; \sS)) \rightarrow M(D; \sS).
$$
Let $f:D \rightarrow D$ be an orientation preserving homeomorphism of the disk
with a single periodic orbit $\sS$. 
The isotopy class of $f$ relative to $\sS$ is represented by
$\beta Z(\mathcal{B}(D; \sS))$ for some braid $\beta = \beta(f,\sS) \in \mathcal{B}(D;\sS)$.
The {\it braid type of $\sS$ for $f$}, denoted by $bt(\sS,f)$ is  the conjugacy class
$\langle \beta Z(\mathcal{B}(D; \sS))\rangle$ in the group $\mathcal{B}(D;\sS)/Z(\mathcal{B}(D; \sS))$.
To simplify the notation of the braid type,
we will write $\langle \beta \rangle$ for $\langle \beta Z(\mathcal{B}(D; \sS))\rangle$.
Let
\begin{eqnarray*}
bt(f) &=& \{bt(P,f)\ |\ P\ \mbox{is\ a\ single\ periodic\ orbit\ of\ }f\},
\end{eqnarray*}
and let $BT$ be the set of all braid types for all homeomorphisms $f:D \rightarrow D$.
A relation $\succeq$ on $BT$ is defined as follows:
for $b_i \in BT$ ($i=1,2$),
\begin{center}
$b_2 \succeq b_1 \Longleftrightarrow$
(For any $f:D \rightarrow D$, $b_2 \in bt(f) \Rightarrow b_1
\in bt(f))$.
\end{center}
We say that $b_2$ {\it forces} $b_1$ if $b_2 \succeq b_1$.
It is known that $\succeq $ gives a partial order on $BT$ \cite{Boy92}, and
we call the relation the {\it forcing relation}.

The topological entropy gives a measure of orbits complexity for a continuous map
of the compact space \cite{Wal}. Let $h(f) \ge 0$ be the topological entropy of $f$.
Let $\beta \in \sB(D;\sS)$ be a pseudo-Anosov braid, and let
$\Phi_{\beta}$ the pseudo-Anosov representative of $\phi_{\beta} \in \sM(D; \sS)$.
Then $\log(\lambda(\Phi_{\beta}))$ is equal to $h(\beta)$, which in turn is the
smallest $h(f)$ among all $f$ with an invariant set $\sS$ such that
$$
bt(\sS,f) = \langle \beta \rangle.
$$
(See \cite{FLP}.)

From the definition of the forcing relation, we have the following.

\begin{cor}
\label{cor_forcing}
If $\beta_1$ and $\beta_2$ are pseudo-Anosov, and  $\langle \beta_2 \rangle \succeq \langle \beta_1 \rangle$,  then
$\lambda(\beta_2) \ge \lambda(\beta_1)$.
\end{cor}

\begin{figure}[htbp]
\begin{center}
\includegraphics[angle=0,width=9cm]{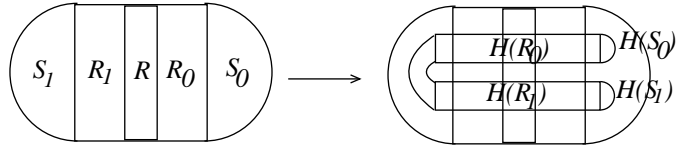}
\caption{Smale horseshoe map.}
\label{fig_smale-h}
\end{center}
\end{figure}

The Smale-horseshoe map $\mathtt{H}: D \rightarrow D$ is a diffeomorphism of the disk
such that the action of $\mathtt{H}$ on
three rectangles $R_0,R_1$, and $R$ and two half disks $S_0,S_1$
is given in Figure \ref{fig_smale-h} \cite{Smale67}.
The restriction of $\mathtt{H}$ to $R_i$ ($i=0,1$) is an affine map such that
$\mathtt{H}$ contracts $R_i$ vertically and stretches horizontally.
The restriction of $\mathtt{H}$ to $S_i$ ($i=0,1$) is a contraction map.
Katok shows that
any $C^{1 + \epsilon}$ surface diffeomorphism ($\epsilon >0$) with  positive topological entropy has
a horseshoe  in some iterate \cite{Katok80}, and
this suggests that the Smale-horseshoe map is a fundamental model for chaotic dynamics.

The set  $\Omega = \displaystyle\bigcap_{n \in {\Z}} \mathtt{H}^n (R_0 \cup R_1)$ is invariant
under $\mathtt{H}$.
Let $\Sigma_2= \{0,1\}^{\Z}$ and let
\begin{eqnarray*}
\sigma: \Sigma_2 &\rightarrow& \Sigma_2
\\
(\cdots w_{-1}\cdot w_0 w_1 \cdots) &\mapsto& (\cdots w_{-1}w_0 \cdot w_1 \cdots),
\quad w_i \in \{0,1\}
\end{eqnarray*}
be the {\it shift map}.
There is a conjugacy $\mathcal{K}: \Omega \rightarrow \Sigma_2$
between the two maps
$\mathtt{H}|_{\Omega}: \Omega \rightarrow \Omega$ and
$\sigma: \Sigma_2 \rightarrow \Sigma_2$ as follows:
\begin{eqnarray*}
\mathcal{K}: \Omega &\rightarrow& \Sigma_2
\\
x &\mapsto& (\cdots \mathcal{K}_{-1}(x) \mathcal{K}_0(x) \mathcal{K}_1(x) \cdots),
\end{eqnarray*}
where
\[
\mathcal{K}_i(x) =
\left\{
\begin{array}{ll}
0 \hspace{3mm}\  \mbox{if\ }&\mathtt{H}^i(x) \in R_0,
\\
1 \hspace{3mm}\  \mbox{if\ }&\mathtt{H}^i(x) \in R_1.
\end{array}
\right.
\]
If $x$ is a period $k$ periodic point, then
the finite word $ (\mathcal{K}_0(x) \mathcal{K}_1(x) \cdots, \mathcal{K}_{k-1}(x))$ is called
{\it code} for $x$.
We say that a braid $\beta $ on the disk  is a {\it horseshoe braid}
if there is a periodic orbit for the Smale-horseshoe map whose braid type is
$\langle \beta \rangle$.
We define a horseshoe braid type in a similar manner.
For the study of  the restricted forcing relation on the set of horseshoe braid types, see
\cite{Hall94},\cite{dCH},\cite{dCH02}, and \cite{dCH03}.

\begin{figure}[htbp]
\begin{center}
\includegraphics[angle=0,width=6cm]{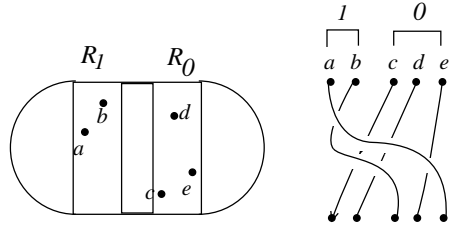}
\caption{Periodic orbit with code $10010$, and its braid representative.}
\label{fig_smale-h-orbit}
\end{center}
\end{figure}

The result by Katok together with Corollary~\ref{cor_forcing} implies that
horseshoe braids are relevant candidates realizing the least dilatation.
We note that  the braid $\sigma_{m,n}$ for each $n \ge m+2$ is a horseshoe braid.
More precisely, the braid type of the periodic orbit with code
$1 \underbrace{0 \cdots 0}_{n-1} 1 \underbrace{0 \cdots 0}_{m}$ or
$1 \underbrace{0 \cdots 0}_{n-1} 1 \underbrace{0 \cdots 0}_{m-1}1$
is represented by the conjugate of $\sigma_{m,n}$ shown in Figure \ref{pA-fig}(c).
Figure~\ref{fig_smale-h-orbit}  illustrates the periodic orbit with code 10010
and its braid representative.
Theorem 15(b) in \cite{dCH03} directly shows the following.
\begin{cor}
If $n \ge n' \ge m+2$, then
$bt(\sigma_{m,n}) \succeq bt(\sigma_{m,n'})$.
\end{cor}

By Propositions \ref{min-prop} and \ref{min-Smn-Bmm-prop},
it is natural to ask whether
the braid $\sigma_{m-1,m+1}$ (respectively,  $\sigma_{m-1,m+2}$)
have the least dilatation among the horseshoe braids on
$(2m+1)$ strands (respectively,  $(2m+2)$ strands).
For the case of even strands,
there is a horseshoe braid having dilatation less than $\lambda(\sigma_{m-1,m+2})$.
In fact, the braid type of period $8$ periodic orbit with code $10010110$ is given by
$$ \langle \beta= \sigma_7 \sigma_6 \sigma_7
 \sigma_5 \sigma_6 \sigma_7
  \sigma_4 \sigma_5 \sigma_6 \sigma_7
 \sigma_3 \sigma_4 \sigma_5 \sigma_6
 \sigma_2 \sigma_3 \sigma_4
\sigma_1 \sigma_2 \rangle,
$$ which satisfies
$\lambda(\beta)= 1.4134< \lambda(\sigma_{2,5}) = 1.5823$, see \cite[Table 2]{Hall94}.
\begin{ques}
For every pseudo-Anosov horseshoe braid $\beta \in \sB(D, 2g+1)$, is it true that
$$
\lambda(\sigma_{g-1,g+1}) \le \lambda(\beta)
\quad ?
$$
\end{ques}

\subsection{Fibered links}\label{fiberedQ-section}

Given a fibered link $(K,F)$, with fibering surface $F$, the {\it monodromy} 
$$\phi_{(K,F)} : F \rightarrow F$$
is the map defined up to isotopy such that the complement of a regular neighborhood of
$K$ in $S^3$ is a mapping torus for $\phi_{(K,F)}$.
Define $\Delta_{(K,F)}$ to be the characteristic
polynomial for the monodromy $\phi_{(K,F)}$ restricted to first homology $\HH_1(F,\R)$.
If $K$ is a fibered knot, then  $\Delta_{(K,F)}$ is the Alexander polynomial of $K$
\cite{Rolfsen76}, \cite{Kaw:Survey}.

For a polynomial $f$,
let $\lambda(f)$ be the maximum norm of roots of $f$ (as in Section~\ref{Salem-Boyd-section}). 
The {\it homological dilatation} of a pseudo-Anosov map $\phi: F \rightarrow F$ is defined to
be $\lambda(f)$, where $f$ is the characteristic polynomial for the restriction of $\phi$ to
$\HH_1(F;\R)$.  Thus, if $(K,F)$ is a fibered link and $\phi_{(K,F)}$ is the
monodromy, then $\lambda(\Delta_{(K,F)})$ is the homological dilatation
of $\phi_{(K,F)}$.  The homological and geometric dilatation are equal
if $\phi$ has orientable stable and unstable foliations \cite{Rykken99}.  We will call
$\phi$ orientable in this case.

Any monic reciprocal integer polynomial is equal to $\Delta_{(K,F)}$ for some
fibered link $(K,F)$ up to multiples of $(t-1)$ and $\pm t$ \cite{Kanenobu81}.
In particular, any Perron polynomial can be realized.    On the other hand, if
the monodromy $\phi_{(K,F)}$ does not have orientable invariant foliations,  then 
$\lambda(\phi_{(K,F)})$ is in general strictly greater than 
 $\lambda(\Delta_{(K,F)})$.

For $g=5$,  Leininger \cite{Leininger04} exhibited a pseudo-Anosov
map $\psi_L$ with dilatation $\lambda_L$ where 
$$
\log (\lambda_L) = 0.162358. 
$$
A comparison shows that this number is strictly less than our candidate for least
element of the braid spectrum $\Sigma (\sB(D,{2g+1})$, for $g=5$:
$$
\log(\lambda(\sigma_{4,6}) )= 0.240965.
$$
The mapping class
$\psi_L$ is realized as the monodromy of the fibered $(-2,3,7)$-pretzel link.
Its dilatation $\lambda_L$ is the smallest known Mahler measure
among monic integer polynomials \cite{Lehmer33}, \cite{Boyd81}.  

In the rest of this section, we will construct a fibered link whose monodromy is obtained
from the braids $\beta_{m,n}$.  
Let $g = \lfloor{\frac{m+n}{2}}\rfloor$, and let 
$$
i= 
\left \{
\begin{array}{ll}
1&\qquad\mbox{if $m+n$ is even},\\
2&\qquad\mbox{if $m+n$ is odd.}
\end{array}
\right .
$$
Let $\sS$ be the set of marked points on $D$ corresponding to the strands of $\beta_{m,n}$,
and let $F$ be the double covering of $D$, branched over $\sS$.  Then $F$ has 
one boundary component if $m+n$ is even and two if $m+n$ is odd.
Let $\phi'_{m,n}$ be the lift of $\phi_{\beta_{m,n}}$
to $F$.  Using an argument argument similar to that in
the proof of Proposition~\ref{spectrum-prop}, we have
$$
\lambda(\phi'_{m,n}) = \lambda(\beta_{m,n}).
$$
Note that  $\phi'_{m,n}$ is one-pronged near
each of the boundary of $F$ if $m+n$ is odd.

\begin{figure}[htbp]
\begin{center}
\includegraphics[height=1.5in]{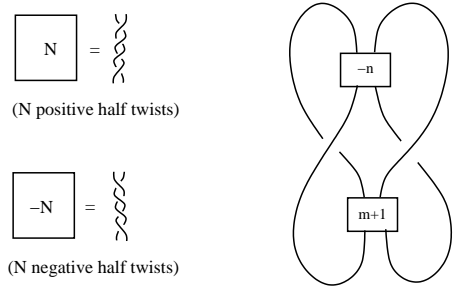}
\caption{Two bridge link associated to $\beta_{m,n}$.}
\label{twobridge-fig}
\end{center}
\end{figure}

Let $K_{m,n}$ be the two-bridge link given in Figure~\ref{twobridge-fig}.
By viewing $(S^3,K_{m,n})$ as the result of a sequence
of Hopf plumbings (see \cite[Section 5]{Hironaka:Salem-Boyd}) one has the following.

\begin{prop}  The complement of a regular neighborhood of $K_{m,n}$ in $S^3$ is
a mapping torus for $\phi'_{m,n}$.
\end{prop}

\noindent
The fibered links $K_{m,n}$ and the dilatations of $\phi'_{m,n}$ were also studied in
\cite{Brinkmann04}.

Let $\Delta_{m,n}$ be the Alexander polynomial for the fibered link $K_{m,n}$.
Salem-Boyd sequences for $\Delta_{m,n}$ were computed in \cite{Hironaka:Salem-Boyd}.
Proposition~\ref{Bmn-lift-prop} implies the following.

\begin{lem} If $m$ and $n$ are both odd, then
$$
\lambda(\beta_{m,n}) = \lambda(\phi'_{m,n}) =  \lambda(\Delta_{K_{m,n}}).
$$
\end{lem}

\begin{ques} Is there a fibered link whose mapping torus is a lift of
$\phi_{\sigma_{m,n}}$?
\end{ques}

\bibliographystyle{math}
\bibliography{math}

\vspace{1in}

\noindent
Eriko Hironaka\newline
Department of Mathematics\newline
Florida State University\newline
Tallahassee, FL 32306-4510\newline
U.S.A.\newline
\bigskip

\noindent
Eiko Kin\newline
Department of Mathematics\newline
 Faculty of Science\newline
 Kyoto University\newline
 Kyoto 606-8502\newline
 Japan\newline

\end{document}